\documentclass[12pt]{article}

\usepackage{amsfonts}
\usepackage{xr}
\usepackage{graphicx}
\usepackage{amsmath}

\def\Z{\mathbb{Z}}
\def\R{\mathbb{R}}
\def\N{\mathbb{N}}

\def\tilde{\widetilde}
\def\epsilon{\varepsilon}
\def\trait (#1) (#2) (#3){\vrule width #1pt height #2pt depth #3pt}
\def\fin{\hfill\trait (0.1) (5) (0) \trait (5) (0.1) (0) \kern-5pt \trait (5) (5) (-4.9) \trait (0.1) (5) (0)}

\newtheorem{thm}{\bf Theorem}[section]
\newtheorem{lem}[thm]{\bf Lemma}
\newtheorem{prop}[thm]{\bf Proposition}

\newtheorem{defi}[thm]{\bf Definition}
\newtheorem{hyp}{\bf Hypothesis}

\title{Traveling fronts in space-time periodic media}
\date{}
\author{Gr\'egoire Nadin\thanks{
D\'epartement de Math\'ematiques et Applications,
\'Ecole Normale Sup\'erieure, CNRS UMR8553 ,
            45 rue d'Ulm, F~75230 Paris cedex 05, France ; e-mail: nadin@dma.ens.fr}}

\begin{document}
\maketitle

\noindent {\bf Abstract:}
This paper is concerned with the existence of pulsating traveling fronts for the equation:
\begin{equation}
\partial_{t} u - \nabla \cdot(A(t,x) \nabla u)+q(t,x)\cdot\nabla u=f(t,x,u),
\end{equation}
where the diffusion matrix $A$, the advection term $q$ and the
reaction term $f$ are periodic in $t$ and $x$. We prove that there exist some speeds $c^{*}$ and $c^{**}$ such that there exists a pulsating traveling front of speed $c$ for all $c\geq c^{**}$ and that there exists no such front of speed $c<c^{*}$. We also give some spreading properties for front-like initial data. In the case of a KPP-type reaction term, we prove that $c^{*}=c^{**}$ and we characterize this speed with the help of a family of eigenvalues associated with the equation. If $f$ is concave with respect to $u$, we prove some Lipschitz continuity for the profile of the pulsating traveling front.

\bigskip

\noindent {\bf R\'esum\'e :}
Cet article \'etudie l'existence de fronts pulsatoires pour l'\'equation :
\begin{equation}
\partial_{t} u - \nabla \cdot(A(t,x) \nabla u)+q(t,x)\cdot\nabla u=f(t,x,u),
\end{equation}
o\`u la matrice de diffusion $A$, le terme d'advection $q$ et le terme de r\'eaction $f$ sont p\'eriodiques en $t$ et en $x$. Nous prouvons l'existence de deux vitesses $c^{*}$ et $c^{**}$ telles qu'il existe un front pulsatoire de vitesse $c$ pour tout $c\geq c^{**}$ et qu'il n'existe pas de tel front de vitesse $c<c^{*}$. Nous donnons \'egalement des propri\'et\'es de $spreading$ pour des donn\'ees initiales ressemblant \`a des fronts. Dans le cas d'un terme de r\'eaction de type KPP, nous prouvons que $c^{*}=c^{**}$ et nous caract\'erisons cette vitesse \`a l'aide d'une famille de valeurs propres associ\'ee \`a l'\'equation. Si $f$ est concave en $u$, nous montrons que le profil du front pulsatoire construit est 
lipschitzien. 

\bigskip


\bigskip


\section{Introduction and preliminaries}

\subsection{Introduction}

This paper investigates the equation:
\begin{equation} \label{eqprinc}
\partial_{t} u - \nabla \cdot(A(t,x) \nabla u)+q(t,x)\cdot\nabla u=f(t,x,u)
\end{equation}
where the coefficients are periodic in $t$ and in $x$. This is a generalization of the homogeneous reaction-diffusion equation:
$ 
\partial_{t} u - \Delta u=u(1-u)
$, which has been first investigated in the pioneering articles of Kolmogorov, Petrovsky and Piskunov \cite{KPP} and Fisher \cite{Fisher}.

The behaviour of the solutions of the homogeneous equation is interesting. First, there exist $planar \ fronts$, that is, solutions of the form $u(t,x)=U(x\cdot e+ct)$, where $e$ is a unit vector and $c$ is the speed of propagation in the direction $-e$. Next, beginning with a positive initial datum $u_{0} \not\equiv 0$ with compact support, we get $u(t,x)\to 1$ when $t\to+\infty$, locally in $x$, moreover, the set where $u$ is close to $1$ spreads with a speed which is equal to the minimal speed of the planar fronts (see \cite{Aronson}) in dimension $1$.

The equation $(\ref{eqprinc})$ arises in population genetics, combustion and population dynamics models. The existence of fronts and the spreading properties have useful interpretations. In population dynamics models, it is very relevant to consider heterogeneous environments and to study the effect of the heterogeneity on the propagation properties.
The homogeneous equation has been fully investigated, but the study of propagation phenomenas for $heterogeneous$ equations is quite recent. It has started with the articles of Freidlin and Gartner \cite{Gartner} and Freidlin \cite{Freidlin2}, who have investigated propagation phenomenas in space periodic environments. They used a stochastic method and avoided the proof of the existence of fronts.

Next, in \cite{Biologicalinvasions, Shigesada1}, Shigesada, Kawasaki and Teramoto defined the notion of $pulsating$ $traveling$ $fronts$, which is a generalization of the notion of planar fronts to space periodic environments. Namely, a solution of equation $(\ref{eqprinc})$, where $A, q$ and $f$ do not depend on $t$ and are $L_i$-periodic with respect to $x_i$, is a pulsating traveling front if $u$ satisfies:
\begin{equation} \label{travelingfrontsespace} \left\{ \begin{array}{l}
\forall x \in \R, t \in \R, u(t+\frac{L\cdot e}{c},x)=u(t,x+L),\\
u(t,x) \rightarrow 0 \ \hbox{as} \ x\cdot e \rightarrow -\infty \ \hbox{and} \ u(t,x) \rightarrow 1 \ \hbox{as} \ x\cdot e \rightarrow +\infty, \\
\end{array} \right. \end{equation}
where $L=(L_1,...,L_N)$, $c$ is the speed of propagation and $0$ and $1$ are the unique zeros of the reaction term $f(x,\cdot)$ for all $x$. 
They did not prove any analytical result but carried out numerical approximations and heuristic computations. One can easily remark that it is equivalent to say that $u$ is a pulsating traveling front if it can be written $u(t,x)=\phi(x\cdot e+ct,x)$, where:
\begin{equation} \label{travelingfrontsespacephi} \left\{ \begin{array}{l}
(z,x)\mapsto\phi(z,x) \ \hbox{is periodic in} \ x,\\
\phi(z,x)\rightarrow 0 \ \hbox{as} \ z \rightarrow -\infty \ \hbox{and} \ \phi(z,x) \rightarrow 1 \ \hbox{as} \ z\rightarrow +\infty. \\
\end{array} \right. \end{equation}

Berestycki and Nirenberg \cite{Nirenberg} and Berestycki and Larrouturou and Lions \cite{Larrouturou} have proved the existence of traveling fronts for heterogeneous advection which does not depend on the variable of the direction of propagation. This result has been generalized, using the notion of pulsating traveling fronts, to the case of a space periodic advection in \cite{Xin} and, next, to the case of a fully space periodic environment with positive nonlinearity in \cite{Frontexcitable, Base2}. It is now being extended to almost periodic environments (see \cite{NadinRossi}).

The existence of fronts has also been proved in the case of a time periodic environment with positive nonlinearity by Fr\'ejacques in \cite{Frejacques}. In this case, the definition of a pulsating traveling front can be easily extended, namely, a solution $u$ is a pulsating traveling front in $t$ if it satisfies:
\begin{equation} \label{travelingfrontstemps} \left\{ \begin{array}{l}
\forall t, x \in \R^{N}, t \in \R, u(t+T,x)=u(t,x+cTe),\\
u(t,x) \rightarrow 0 \ \hbox{as} \ x\cdot e \rightarrow -\infty \ \hbox{and} \ u(t,x) \rightarrow 1 \ \hbox{as} \ x\cdot e \rightarrow -\infty, \\
\end{array} \right. \end{equation}
where $T$ is the period of the coefficients.
The existence of time-periodic pulsating traveling front has also been proved in the case of a bistable nonlinearity in \cite{Alikakos}. In the case of an almost periodic environment with bistable nonlinearity, Shen has defined a pulsating traveling front in \cite{Shen1,Shen2}, the speed of propagation is then almost periodic and not constant.

Let us mention the most recent breakthroughs on this topic to conclude this introduction. Two definitions for a notion of fronts in general media have been given by Berestycki and Hamel \cite{Generalwaves2, Generalwaves} and by Matano \cite{Matanodef}. Using Matano's definition, Shen has proved the existence of such fronts in time general media with bistable nonlinearity \cite{Shengeneral}. Using Berestycki and Hamel's definition, it has been proved in a parallel way by Nolen and Ryzhik \cite{NolenRyzhik} and Mellet and Roquejoffre \cite{MelletRoquejoffre} that such fronts exist in space general media with ignition type nonlinearity.


\subsection{Notion of fronts in space-time periodic media}

The investigation of space-time periodic reaction-diffusion equation is very recent. In 2002, Weinberger has proved the existence of pulsating traveling fronts in  this case in a discrete context \cite{Weinberger}. In 2006, Nolen, Rudd and Xin investigated the case of an incompressible periodic drift in \cite{Nolen, Nolenshear}, with a positive homogeneous nonlinearity $f(u)$. 

In order to define a notion of front in space-time periodic media, one can first try to extend Definition $\ref{travelingfrontsespace}$ and say that a pulsating traveling front will be a solution of $(\ref{eqprinc})$ which satisfies some equality
\[u(t+\frac{pL}{c},x)=u(t,x+pL)\]
for all $(t,x)\in\R\times\R^{N}$ and for some $p\in\Z$ since $pL$ is a space period of the medium. But it has been proved in \cite{Nolen} that this necessarily implies the existence of some $q\in\Z\backslash \{0\}$ such that $c=\frac{pL}{qT}$, which is not satisfactory since we expect to find a half-line of speeds associated with fronts and not only a sequence, like in space periodic or time periodic media.

One can then try to extend definition $(\ref{travelingfrontsespacephi})$ and say that a front is a solution of $(\ref{eqprinc})$ which can be written $u(t,x)=\phi(x\cdot e+ct,t,x)$, where $(z,t,x)\mapsto \phi(z,t,x)$ is periodic in $t$ and $x$ and converges to $1$ as $z\rightarrow+\infty$ and to $0$ as $z\rightarrow-\infty$. Such a function $\phi$ has to satisfy:
\begin{equation} \label{eqtf} \begin{array}{l}\partial_{t} \phi - \nabla \cdot (A(t,x) \nabla \phi) - eA(t,x)e \partial_{zz}\phi -2e A(t,x)\partial_{z}\nabla \phi\\
+q(t,x)\cdot\nabla\phi+q(t,x)\cdot e\partial_{z}\phi+c\partial_{z}\phi=f(t,x,\phi),\\ \end{array} \end{equation}
over the hyperplane $z=x\cdot e+ct$. Thus, in order to prove the existence of such fronts, one can try to find a solution of this equation over the whole space $\R\times\R\times\R^{N}$ and to set $u(t,x)=\phi(x\cdot e+ct,t,x)$. But then some difficulties arise.

Actually, equation  $(\ref{eqtf})$ is degenerate and thus it is not easy to construct a regular solution of (\ref{eqtf}). If one defines $v(y,t,x)=\phi (x\cdot e+ct+y,t,x)$, then $\phi$ has the same regularity as $v$ and this function satisfies:
\[\partial_{t}v-\nabla\cdot (A(t,x)\nabla v)+q(t,x)\cdot \nabla v=f(t,x,v) \ \hbox{in} \ \mathcal{D}'(\R\times\R\times\R^{N}).\]
As this equation does not depend on $y$, one cannot expect to get some regularity in $y$ from it. For example the Hormander-Kohn conditions (see \cite{Hormander, Kohn}) do not hold here. Thus $\phi$ may only be measurable and not continuous. In this case, setting $u(t,x)=\phi (x\cdot e+ct,t,x)$ is not relevant since the hyperplane $z=x\cdot e+ct$ is of measure zero and thus many functions $u$ may be written in this form. 

We underline that this kind of issue does not arise in space periodic or time periodic media since in these cases, we do not add an extra-variable $z$. Thus we can go from $\phi$ to $u$, which satisfies the regularizing equation $(\ref{eqprinc})$, and then go back to $\phi$ which then has the same regularity as $u$. 

\bigskip

Because of this issue, Nolen, Rudd and Xin gave a weakened definition of pulsating traveling fronts. This definition was stated in \cite{Nolen} in the case of a space-time periodic incompressible advection, with homogeneous $A$ and $f$, but it can be naturally extended to equation $(\ref{eqprinc})$:
\begin{defi} \cite{Nolen}\label{deftravelingfrontsNolen} 
Assume that equation $(\ref{eqprinc})$ admits two space-time periodic solutions $p^{-}$ and $p^{+}$ such that $p^{-}(t,x) <p^{+}(t,x)$ for all $(t,x)\in\R\times\R^{N}$.
Then a front traveling at speed $c$ is a function $\phi(z,t,x)\in L^{\infty}_{loc}(\R\times\R\times\R^{N})$ whose directional
derivatives $\partial_{t}\phi+c\partial_{z}\phi$, $\nabla\phi+e\partial_{z}\phi$ and $(\nabla+e\partial_{z})^{2}\phi$ are continuous and satisfy the equation:
\begin{equation} \label{eqNolen} \begin{array}{l}\partial_{t} \phi - \nabla \cdot (A(t,x) \nabla \phi) - eA(t,x)e \partial_{zz}\phi - \nabla \cdot (A(t,x)e \partial_{z} \phi)\\
-\partial_{z}(e A(t,x)\nabla \phi)+q(t,x)\cdot\nabla\phi+q(t,x)\cdot e\partial_{z}\phi+c\partial_{z}\phi=f(t,x,\phi) \ \hbox{in} \ \mathcal{D}'(\R\times\R\times\R^{N}),\\ \end{array} \end{equation}
such that $\phi$ is periodic in $t$ and $x$ and 
\begin{equation} \left\{ \begin{array}{l}
\phi(z,t,x)-p^{-}(t,x) \rightarrow 0 \ \hbox{as} \ z\rightarrow -\infty \ \hbox{uniformly in} \ (t,x) \in \R\times\R^{N},\\
\phi(z,t,x)-p^{+}(t,x) \rightarrow 0 \ \hbox{as} \ z\rightarrow +\infty \ \hbox{uniformly in} \ (t,x) \in \R\times\R^{N}.\\
\end{array} \right. \end{equation}\end{defi}

The difficulty with this definition is that it is not clear if it has a direct link with the parabolic equation $(\ref{eqprinc})$ since the profile $\phi$ is only measurable and thus $u(t,x)=\phi(x\cdot e+ct,t,x)$ does not really make sense. This is why we now give the following equivalent definition:
\begin{defi} \label{deftravelingfronts}
We say that a function $u$ is a $pulsating \ traveling \ front$ of speed $c$ in the direction $-e$ that connects $p^{-}$ to $p^{+}$ if it can be written $u(t,x)= \phi(x\cdot e+ct,t,x)$, where $\phi\in L^{\infty}(\R\times\R\times\R^{N})$ is such that for almost every $y\in\R$, the function $(t,x)\mapsto\phi(y+x\cdot e+ct,t,x)$ satisfies equation $(\ref{eqprinc})$. We ask the function $\phi$ to be periodic in its second and third variables and to satisfy:
\begin{equation} \left\{ \begin{array}{l}
\phi(z,t,x)-p^{-}(t,x) \rightarrow 0 \ \hbox{as} \ z\rightarrow -\infty \ \hbox{uniformly in} \ (t,x) \in \R\times\R^{N},\\
\phi(z,t,x)-p^{+}(t,x) \rightarrow 0 \ \hbox{as} \ z\rightarrow +\infty \ \hbox{uniformly in} \ (t,x) \in \R\times\R^{N}.\\
\end{array} \right. \end{equation}
\end{defi}

\noindent {\bf Remark:} The equivalence between the two definitions is not obvious and will be proved later.

\bigskip

Of course these definitions are not very convenient and we would like to construct pulsating traveling fronts that are at least continuous. We will prove in this article that this is possible under some KPP-type assumption. 
\begin{defi} \label{defcontptf}
We say that a solution $u$ of $(\ref{eqprinc})$ is a $Lipschitz \ continuous$ pulsating traveling front of speed $c$ in the direction $-e$ that connects $p^{-}$ to $p^{+}$ if it can be written $u(t,x)= \phi(x\cdot e+ct,t,x)$, where $\phi\in W^{1,\infty}(\R\times\R\times\R^{N})$ is periodic in its second and third variables and satisfies:
\begin{equation} \left\{ \begin{array}{l}
\phi(z,t,x)-p^{-}(t,x) \rightarrow 0 \ \hbox{as} \ z\rightarrow -\infty \ \hbox{uniformly in} \ (t,x) \in \R\times\R^{N},\\
\phi(z,t,x)-p^{+}(t,x) \rightarrow 0 \ \hbox{as} \ z\rightarrow +\infty \ \hbox{uniformly in} \ (t,x) \in \R\times\R^{N}.\\
\end{array} \right. \end{equation}
\end{defi}

One can check that Definition \ref{defcontptf} fits with the definition of a $generalized$ $almost$ $planar$ $traveling$ $wave$ of speed $c$ that has been given by H. Berestycki and F. Hamel in \cite{Generalwaves2, Generalwaves} with, using the notations of this reference, $\Gamma(t)=\{x\in\R^{N}, x\cdot e+ct=0\}$.


\subsection{Framework}

In \cite{Nolen}, Nolen, Rudd and Xin proved the following results:
\begin{thm} \label{thmNolen} \cite{Nolen}
Assume that $f$ and $A$ do not depend on $(t,x)$, that $\nabla\cdot q\equiv 0$, that $\int_{(0,T)\times C}q=0$ and that $f(0)=f(1)=0$, $f'(0)>0$, $f'(1)<0$ and for all $s\in (0,1)$, $f(s)>0$. Then:

1) there exists a speed $c^{*}$ such that there exists a pulsating traveling front of speed $c^{*}$,

2) if $f$ is of KPP type, that is, $f(s)\leq f'(0)s$ for all $s>0$, then there exists no pulsating traveling fronts of speed $c<c^{*}$. Furthermore, the speed $c^{*}$ can be characterized with the help of some space-time periodic principal eigenvalues associated with the problem.
\end{thm}

This theorem leaves some open questions. First of all, does it exist a pulsating traveling front $for \ all$ $c>c^{*}$, which is the classical result for monostable nonlinearities in space periodic or time periodic media (see \cite{Frontexcitable, Base2, Frejacques})? Secondly, is this possible to construct regular profiles that are associated with solutions of the parabolic equation (\ref{eqprinc})? And lastly, one can wonder if these results can be extended to the case of a heterogeneous diffusion matrix $A$ and reaction term $f$. 

\bigskip

In the present paper, we prove the existence of such pulsating traveling fronts when not only $q$, but also $A$ and $f$ are periodic in $t$ and $x$, with different methods as that of \cite{Nolen}. These new methods enable us to prove that there exists a half-line of speeds $[c^{**},+\infty)$ associated with pulsating traveling fronts, which is a stronger result than the result of \cite{Nolen}. We also prove that there exists some speed $c^{*}$ associated with some pulsating traveling front such that there exists no pulsating traveling front of speed $c<c^{*}$, even if $f$ is not of KPP type. If $f$ satisfies a KPP type assumption, then $c^{*}=c^{**}$ and we characterize the minimal speed $c^{*}$. 

In addition, we investigate new questions in this article. We prove in particular that if $s\mapsto f(t,x,s)/s$ is nonincreasing, then there exist some \textit{Lipschitz continuous pulsating traveling fronts} of speed $c$ if and only if $c\geq c^{*}$. 

Except in \cite{Base2, Weinberger}, all the preceding papers, including \cite{Nolen}, only considered the case of positive nonlinearities. This case makes sense in the case of a combustion model, but not in populations dynamics models. Usually, in this kind of models, the reaction term has the form $f(t,x,u)=u(\mu(t,x) -u)$, where $\mu$ is the difference between a birth rate and a death rate when the population is small, which both depend on the environment. In unfavourable areas, this term may be negative.
Moreover, such a reaction term often leads to heterogeneous asymptotic states. In the present paper, we do not necessarily assume that the reaction term $f$ is positive.

Lastly, all the previous papers were only considering incompressible drifts of null-average. As such drifts do not exist in dimension $1$, it was previously impossible to study the effect of the advection on the propagation of the fronts in dimension $1$. In the sequel, we will not make any such assumption on the drift term $q$. The dependence between the spreading speed and a compressible advection term has been investigated by Nolen and Xin in \cite{NolenXin1d} and the author in \cite{dependance}, who proved that such an advection term may decrease the spreading speed. 

To sum up, we are going to prove the existence of pulsating traveling fronts in space-time periodic environments under very weak hypotheses that only rely on the stability of the steady state $0$ and on the uniqueness of the space-time periodic asymptotic state $p$. Hence, this paper gives new results even in space periodic or time periodic environments.

These results are summarized in the note \cite{notetf}.

\subsection{Hypotheses}

We will need some regularity assumptions on $f, A, q$.
The function $f: \mathbb{R} \times \mathbb{R}^{N} \times
\mathbb{R}^{+} \rightarrow \mathbb{R}$ is supposed to be of class $C^{\frac{\delta}{2},\delta}$ in $(t,x)$ locally in $u$ for a given $0<\delta<1$ and of class
$C^{1,r}$ in $u$ on $\mathbb{R}\times\mathbb{R}^{N}\times [0,\beta]$ for some given $\beta>0$ and $r>0$. We also assume that $0$ is a state of equilibrium, that is, $\forall x,\forall t,  f(t,x,0)=0$.

The matrix field $A:\mathbb{R} \times \mathbb{R}^{N}
\rightarrow S_{N}(\mathbb{R})$ is supposed to be of class $C^{\frac{\delta}{2},1+\delta}$. We suppose furthermore that $A$ is uniformly elliptic and continuous: there exist some positive constants $\gamma$ and $\Gamma$ such that for all $\xi \in \mathbb{R}^{N}, (t,x) \in \mathbb{R} \times \mathbb{R}^{N}$ one has:
\begin{equation} \label{ellipticity} \gamma \|\xi\|^{2} \leq \sum_{1 \leq i,j \leq N} a_{i,j}(t,x) \xi_{i} \xi_{j} \leq \Gamma \|\xi\|^{2}, \end{equation}
where $\|\xi\|^{2}=\xi_{1}^{2}+...+\xi_{N}^{2}$.

The drift term $q: \mathbb{R} \times \mathbb{R}^{N} \rightarrow \mathbb{R}^{N}$ is supposed to be of class $C^{\frac{\delta}{2},\delta}$ and we assume that $\nabla\cdot q\in L^{\infty}(\R\times\R^{N})$.

Moreover, we assume that $f$, $A$ and $q$ are periodic in $t$ and $x$. That is, there exist some positive constant $T$ and some vectors $L_{1},...,L_{N}$, where $L_{i}$ is colinear to the axis of coordinates $e_{i}$, for all $(t,x,s)\in\R\times\R^{N}\times\R^+$ and for all $i\in [1,N]$, one has:
\[\begin{array}{rccccccccc}
A(t+T,x)&=&A(t,x),\ q(t+T,x)&=&q(t,x) &\hbox{ and }& f(t+T,x,s)&=&f(t,x,s),\\
A(t,x+L_i)&=&A(t,x), \ q(t,x+L_i)&=&q(t,x) &\hbox{ and }& f(t,x+L_i,s)&=&f(t,x,s).\\
\end{array}\]
We define the periodicity cell $C= \Pi_{i=1}^{N} (0,|L_{i}|)$.

\bigskip

The only strong hypothesis that we need is the following one:
\begin{hyp}\label{uniquenesshyp}
Equation $(\ref{eqprinc})$ admits a positive continuous space-time periodic solution $p$. Furthermore, if $u$ is a space periodic solution of equation $(\ref{eqprinc})$ such that $u\leq p$ and $\inf_{(t,x)\in\R\times\R^{N}}u(t,x)>0$, then $u\equiv p$.
\end{hyp}

It is not easy to check that this hypothesis is true. This condition is investigated in section \ref{example-section}. If the solution $p$ is the unique space periodic uniformly positive solution of $(\ref{eqprinc})$, then this hypothesis is satisfied. But we do not need a general uniqueness hypothesis. 

We now define the $generalized \ principal \ eigenvalue$ associated with equation $(\ref{eqprinc})$ in the neighborhood of the steady states $0$:
\begin{equation} \label{geneigenpetite}\begin{array}{ll}\lambda_{1}' = \inf &\{ \lambda \in \mathbb{R}, \ \exists \phi \in \mathcal{C}^{1,2}(\mathbb{R}\times\mathbb{R}^{N})\cap W^{1,\infty}(\R\times\R^{N}), \ \phi > 0, \ \phi \ \hbox{is T-periodic}, \\
&(-\mathcal{L}+\lambda)\phi \geq 0 \ \hbox{in} \ \mathbb{R}\times\mathbb{R}^{N}\}\\
\end{array} \end{equation}
where $\mathcal{L}$ is the linearized operator associated with equation $(\ref{eqprinc})$ in the neighborhood of $0$:
\[\mathcal{L}\phi=\partial_{t}\phi-\nabla\cdot (A(t,x)\nabla\phi)+q(t,x)\cdot \nabla\phi-f_u'(t,x,0)\phi.\]
The properties of this eigenvalue and its link with equation $(\ref{eqprinc})$ has been investigated in \cite{existence, eigenvalue}. We will assume that $\lambda_{1}'<0$, that is, the steady state $0$ is linearly unstable. If $f$ does not depend on $t$ and $x$, then this hypothesis is equivalent to $f'(0)>0$.

The hypothesis $\lambda_{1}'<0$ is optimal if we want to state a result for general $f$, otherwise, if $\lambda_{1}'\geq 0$, Hypothesis \ref{uniquenesshyp} may be contradicted. In \cite{existence}, we proved the following theorem:
\begin{thm} \cite{existence}
 If $\lambda_{1}'\geq 0$ and if for all $(t,x)\in \R\times \R^N$, the growth rate $s\mapsto f(t,x,s)/s$ is decreasing, then there is no nonnegative bounded continuous entire solution of equation $(\ref{eqprinc})$ except $0$.
\end{thm}

In the previous papers that were considering heterogeneous reaction terms, like \cite{Frontexcitable, Base2}, the authors used some assumption in the neighborhood of $p$ like: $s\mapsto f(t,x,s)$ is decreasing in the neighborhood of $p$. In this article, we managed to get rid of this hypothesis.

\bigskip

Lastly, let us underline that, as all these hypotheses are related to local properties of equation $(\ref{eqprinc})$ we can consider other kinds of equations. For example, the next results are true for the reaction-diffusion equation associated with some stochastic differential equation:
\begin{equation} \label{stochform} \partial_{t}u - \alpha_{ij}(t,x)\partial_{ij}u+\beta_{i}(t,x)\partial_{i}u=f(t,x,u),\end{equation}
where $(\alpha_{ij})_{ij}$ is an elliptic matrix field, that is, it satisfies $(\ref{ellipticity})$ and $(\beta_{i})_{i}$ is general vector field. Setting $q_{j}(t,x)=\beta_{j}(t,x)-\partial_{i}\alpha_{ij}(t,x)$ and $A=\alpha$, one is back to equation $(\ref{eqprinc})$. This change of variables was impossible with the hypotheses of the previous papers since it does not necessarily give a divergence-free vector field $q$. Similarly, one can consider 
\begin{equation} \label{divdrift} \partial_{t}u - \nabla \cdot (A(t,x)\nabla u)+\nabla \cdot (q(t,x) u)=g(t,x,u),\end{equation}
by doing the change of variables $f(t,x,s)=g(t,x,s)-(\nabla \cdot q) (t,x)s$.

\subsection{Examples} \label{example-section}

Hypothesis $\ref{uniquenesshyp}$ is far from being easy to check and we now give two classical examples for which this uniqueness hypothesis holds. Our first example is related to biological models and has been investigated in details in a previous article:

\begin{thm} \cite{existence}\label{reactionhypbio}
 If $\lambda_{1}'<0$ and if $f$ satisfies:
\begin{equation} \label{hyp1}
\forall (t,x) \in \R\times\mathbb{R}^{N}, s\mapsto \frac{f(t,x,s)}{s} \ \hbox{is decreasing},
\end{equation}
\begin{equation} \label{hyp2}
\exists M>0, \forall x \in \mathbb{R}^{N}, \forall t \in
\mathbb{R}, \forall s \geq M, f(t,x,s) \leq 0, \end{equation}
then Hypothesis $\ref{uniquenesshyp}$ is satisfied.
\end{thm}
The two hypotheses on the reaction term $f$ both have a biological meaning (see \cite{Base1, Base2}). The first hypothesis means that the intrinsic growth rate $s\mapsto f(t,x,s)/s$ is decreasing when the population
density is increasing. This is due to the intraspecific
competition for resources. The second hypothesis means that
there is a saturation density: when the population is very
important, the death rate is higher than the birth rate and the
population decreases uniformly. We remark that hypothesis $(\ref{hyp1})$ implies that $f$ is of KPP type, that is, for all $(t,x,s)\in\R\times\R^{N}\times\R^{+}$, one has $f(t,x,s)\leq f_u'(t,x,0)s$. The reader will find more precise existence and uniqueness results for reaction-diffusion equation in space-time periodic media under the hypotheses of Theorem $\ref{reactionhypbio}$ in \cite{noteexistence, existence}.

\bigskip

Our second example is related to combustion models. 
\begin{prop} \label{reactionhyppos}
Assume that:

i) $f(t,x,1)=0$ for all $(t,x)\in\R\times\R^{N}$, 

ii) for all $(t,x)\in\R\times\R^{N}$, if $s\in(0,1)$ then one has $f(t,x,s)>0$.

Then $p\equiv 1$ is the only entire bounded solution of equation $(\ref{eqprinc})$ such that 
\[0<\inf_{\R\times\R^{N}}p\leq \sup_{\R\times\R^{N}}p\leq 1.\] 
Thus Hypothesis $\ref{uniquenesshyp}$ is satisfied.
\end{prop}

\noindent {\bf Proof.}
Assume that $u$ is some uniformly positive continuous entire solution of equation $(\ref{eqprinc})$. Set $m=\inf_{\R\times\R^{N}}u>0$ and consider a sequence $(t_{n},x_{n})\in\R\times\R^{N}$ such that $u(t_{n},x_{n})\rightarrow m$. For all $n$, there exist some $(s_{n},y_{n})\in [0,T]\times\overline{C}$ such that $t_{n}-s_{n}\in T\Z$ and $x_{n}-y_{n}\in \Pi_{i=1}^{N}L_{i}\Z$. Up to extraction, one can assume that $s_{n}\rightarrow s_{\infty}$ and $y_{n}\rightarrow y_{\infty}$. Set $u_{n}(t,x)=u(t+t_{n},x+x_{n})$. This function satisfies:
\[\partial_{t}u_{n}-\nabla \cdot (A(t+s_{n},x+y_{n})\nabla u_{n})+q(t+s_{n},x+y_{n})\cdot \nabla u_{n}=f(t+s_{n},x+y_{n},u_{n}).\]
The Schauder parabolic estimates and the periodicity of the coefficients yield that one can extract a subsequence that converges to some function $u_{\infty}$ in $\mathcal{C}^{1,2}_{loc}(\R\times\R^{N})$ that satisfies:
\[\partial_{t}u_{\infty}-\nabla \cdot (A(t+s_{\infty},x+y_{\infty})\nabla u_{\infty})+q(t+s_{\infty},x+y_{\infty})\cdot \nabla u_{\infty}=f(t+s_{\infty},x+y_{\infty},u_{\infty}).\]
Furthermore, one has $u_{\infty}(0,0)=m$ and $u_{\infty}\geq m$. As $f$ is nonnegative in $\R\times\R^N\times [0,1]$, the strong parabolic maximum principle and the periodicity give $u_{\infty}\equiv m$. If $u\not\equiv 1$, then $m<1$, which would contradict the previous equation since $f(t,x,m)>0$. Thus $m\geq 1$ and $u\equiv 1$.
$\Box$

\bigskip

These two examples prove that our hypotheses are very weak and include the classical hypotheses that were usually used in space periodic or time periodic media. It enables us to tconsider more general coefficients. For example, we can consider very general drifts, even $compressible$ ones which is totally new: the only previous papers that where considering compressible drifts were dealing with spreading properties or qualitative properties of the traveling fronts (see \cite{Hamelmonoticity, NolenXin1d}). We can also consider oscillating reaction terms, that admit several ordered steady states. These hypotheses are hardly optimal and the only open problem that remains is the case where there are two space-time periodic solutions of equation $(\ref{eqprinc})$ that cross each other. We underline that the existence of pulsating traveling fronts in space periodic media for general $f$ with $q\equiv 0$ has been proved at the same time, with the same kind of method, by Guo and Hamel in \cite{GuoHamel2}.

One can also notice that it is easily possible to treat the case where the homogeneous solution $0$ is replaced by some space-time periodic solution $p^{-}$. In this case, the uniqueness Hypothesis $\ref{uniquenesshyp}$ is replaced by the following one:
\begin{hyp}\label{uniquenesshyp2}
 Equation $(\ref{eqprinc})$ admits two positive continuous space-time periodic solutions $p^{-}$ and $p^{+}$. Furthermore, if $u$ is a space periodic solution of equation $(\ref{eqprinc})$ such that $u\leq p^{+}$ and $\inf_{\R\times\R^{N}}(u-p^{-})>0$, then $u\equiv p^{+}$.
\end{hyp}
Setting $g(t,x,s)=f(t,x,s+p^{-}(t,x))-f(t,x,s)$, one can easily get back to the case $p^{-}\equiv 0$ and all the results of this paper can easily be generalized using similar changes of variables.

\subsection{The associated eigenvalue problem}

This section deals with the eigenvalues of the operator:
\begin{equation} \label{Llambda} L_{\lambda} \psi = \partial_{t}\psi - \nabla\cdot (A\nabla \psi) -2\lambda A \nabla \psi+q\cdot\nabla\psi-(\lambda A\lambda+\nabla\cdot(A\lambda)+\mu-q\cdot\lambda)\psi,\end{equation}
where $\lambda \in \mathbb{R}^{N}$ and $\psi \in \mathcal{C}^{1,2}(\mathbb{R} \times \mathbb{R}^{N})$. We assume that $A$ and $q$ satisfy the same hypotheses as in the previous part and $\mu \in \mathcal{C}^{\delta/2, \delta}(\R\times\R^{N})$ is a space-time periodic function. The example we will keep in mind is $\mu(t,x)=f_u'(t,x,0)$. 

\begin{defi}
A $space-time \ periodic \ principal \ eigenfunction$ of the operator $L_{\lambda}$ is a function $\psi \in \mathcal{C}^{1,2}(\mathbb{R} \times \mathbb{R}^{N}) $ such that there exists a constant $k$ so that:
\begin{equation} \label{eigen} \left\{\begin{array}{rcl}
L_{\lambda} \psi &=& k \psi \ \hbox{in} \ \R\times\R^N,\\
\psi &>& 0 \ \hbox{in} \ \R\times\R^N,\\
\psi \ &\hbox{is}& T-\hbox{periodic}\\
\psi \ &\hbox{is}& L_i-\hbox{periodic for} \ i=1,...,N.\\
\end{array}\right.
\end{equation}
Such a $k$ is called a space-time periodic principal eigenvalue.
\end{defi}

This family of eigenvalues has been widely investigated in \cite{eigenvalue}.
The following theorem states the existence and the uniqueness of the eigenelements:

\begin{thm}\cite{eigenvalue}
For all $\mu, A, q, \lambda$, there exists a couple $(k,\psi)$ that satisfies (\ref{eigen}). Furthermore, $k$ is unique and $\psi$ is unique up to multiplication by a positive constant. 
\end{thm}

We define $k_{\lambda}(A,q,\mu)= k$ the space-time periodic principal eigenvalue associated with $L_{\lambda}$. The generalized principal eigenvalue is characterized by:

\begin{prop} \cite{eigenvalue}
 One has: $\lambda_{1}'=k_{0}$.
\end{prop}

For all $A,q,\mu$ and $e\in\mathbb{S}^{N-1}$, we define:
\[ c^{*}_{e}(A,q,\mu)=\min \{ c\in \R, \ \hbox{there exists} \ \lambda>0 \
\hbox{such that} \ k_{\lambda e}(A,q,\mu)+\lambda c=0 \}.\]

We will denote $c^{*}(\mu)=c^{*}_{e}(A,q,\mu)$ in the sequel when there is no ambiguity. This quantity arises when one is searching for exponentially decreasing solutions of $(\ref{eqprinc})$. We will see in the sequel that if $f$ is of KPP-type, the minimal speed of the pulsating traveling fronts in direction $-e$ equals $c^{*}_{e}(A,q,\mu)$, where $\mu(t,x)=f_u'(t,x,0)$. If $f$ is not of KPP-type, this is not true anymore, but one can get some estimates for the minimal speed of propagation using $c^{*}_e(A,q,\mu)$ and $c^{*}_e(A,q,\eta)$, where $\eta(t,x)=\sup_{0<s<p(t,x)}f(t,x,s)/s$.

\section{Statement of the main results}

\subsection{Existence of KPP traveling fronts}

Our first result holds in the KPP case:

\begin{thm} \label{thmexistence}
1) Assume that $\lambda_{1}'< 0$, that Hypothesis $\ref{uniquenesshyp}$ is satisfied and that 
\[f(t,x,s)\leq \mu(t,x)s \ \hbox{for all} \ (t,x,s)\in\R\times\R^{N}\times\R^{+},\]
where $\mu(t,x)=f_u'(t,x,0)$. Then for all unit vector $e$, there exists a minimal speed $c^{*}_{e}$ such that for all speed $c\geq c^{*}_{e}$, there exists a pulsating traveling front $u$ of speed $c$ in direction $-e$ that links $0$ to $p$. 
This speed can be characterized:
\begin{equation} \label{characterizationminspeed} c^{*}_{e}=c^{*}_{e}(A,q,\mu)=\min \{ c\in \R, \ \hbox{there exists} \ \lambda>0 \
\hbox{such that} \ k_{\lambda e}(A,q,\mu)+\lambda c=0 \}. \end{equation}
\smallskip
2) Moreover, for all $c\geq c^{*}_{e}$, the profile $\phi$ of the pulsating traveling front $u$ of speed $c$ we construct is nondecreasing almost everywhere in $z$, that is, for almost every $(z_{1},z_{2})\in\R^{2}$ such that $z_{1}\geq z_{2}$, one has $\phi (z_{1},t,x)\geq \phi (z_{2},t,x)$ for all $(t,x)\in\R\times\R^{N}$. Lastly, for all $c>c^{*}_{e}$, the function $\phi$ satisfies:
\begin{equation}\label{expdec} \phi(z,t,x)\sim \psi_{\lambda_{c}}(t,x)e^{\lambda_{c}(\mu) z} \ \hbox{as} \ z\rightarrow-\infty,\end{equation}
uniformly in $(t,x)\in\R\times\R^{N}$, where $\psi_{\lambda_{c}(\mu)}$ is some principal eigenfunction associated with $k_{\lambda_{c}(\mu)e}(A,q,\mu)$.
\end{thm}

\noindent {\bf Remark:}
The quantity $\lambda_{c}(\mu)$ will be defined in Definition $\ref{deflambdac}$. It is roughly the smallest $\lambda>0$ such that $k_{\lambda e}(A,q,\mu)+\lambda c= 0$.

\bigskip

In the other hand, the following proposition gives a lower bound for the speeds which are associated with pulsating traveling fronts:
\begin{prop} \label{thmnonexistence}
 If $\lambda_{1}'\leq 0$ and Hypothesis $\ref{uniquenesshyp}$ is satisfied, then for all $c<c^{*}_{e}(A,q,\mu)$, where $\mu(t,x)=f_u'(t,x,0)$, there exists no pulsating traveling front of speed $c$.
\end{prop}
This proposition is true even if $f$ is not of KPP type. It is in fact the corollary of some spreading properties for front-like initial data. These spreading properties will be stated later. 

\bigskip

The existence of pulsating traveling fronts with prescribed exponential behaviour has been obtained by Bag\`es in \cite{TheseBages} in space periodic media. The method used by Bag\`es is close to our method and he also managed to construct a profile of speed $c=c^*_e$ which satisfies $\phi(z,t,x)\sim -z \psi_{\lambda_{c}}(t,x)e^{\lambda_{c}(\mu) z}$ as $z\rightarrow -\infty$. We did not manage to consider the critical case $c=c^*_e$ here because of technical issues that we will emphasize later. 

It has recently been proved by Hamel in \cite{Hamelmonotonicity} that in space periodic media, $any$ pulsating traveling front of speed $c\geq c^{*}_{e}$ satisfies $(\ref{expdec})$. In the present paper, we only prove that there exist $some$ pulsating traveling fronts which satisfy this property if $f$ is of KPP type.

\subsection{Existence of traveling fronts for general reaction terms}

If $f$ is not of KPP type, then our results are a little bit weaker:

\begin{thm} \label{thmexistencepos}
Assume that $\lambda_{1}'< 0$ and that Hypothesis $\ref{uniquenesshyp}$ is satisfied. Then for all unit vector $e$, for all speed $c\geq c^{*}_{e}(A,q,\eta)$, where $\eta(t,x)=\sup_{0<s<p(t,x)}\frac{f(t,x,s)}{s}$, there exists a pulsating traveling front $u$ of speed $c$ in direction $-e$ that links $0$ to $p$. 
Moreover, the profile $\phi$ is nondecreasing almost everywhere in $z$, that is, for almost every $(z_{1},z_{2})\in\R^{2}$ such that $z_{1}\geq z_{2}$, one has $\phi (z_{1},t,x)\geq \phi (z_{2},t,x)$ for all $(t,x)\in\R\times\R^{N}$.

If $c<c^{*}_{e}(A,q,\mu)$, where $\mu(t,x)=f_u'(t,x,0)$, then there exists no pulsating traveling fronts of speed $c$.
\end{thm}

The non-existence result is a direct consequence of Proposition $\ref{thmnonexistence}$. Thus, in order to get the result that is announced in the abstract, one only sets $c^{*}=c^{*}_e(A,q,\mu)$ and $c^{**}=c^{*}_e(A,q,\eta)$.

In this case we do not manage to prove estimate $(\ref{expdec})$ but it should be underlined that such an estimate holds in space periodic media (see \cite{Hamelmonoticity}) or when the heterogeneity does not depend on the direction of propagation (see \cite{Nirenberg}).

One can wonder where does this $\eta$ come from. In fact, if $f$ is not of KPP type, one can define:
\[g(t,x,s)=s \sup_{r\geq s}\frac{f(t,x,r)}{r}.\]
One has $g\geq f$, $s\mapsto g(t,x,s)/s$ is nonincreasing for all $(t,x)\in\R\times\R^{N}$ and $g_{u}'(t,x,0)=\eta(t,x)$. Thus $g$ is somehow the lowest KPP nonlinearity which lies above $f$. As it is a KPP nonlinearity, its associated minimal speed is $c^{*}_{e}(g)= c^{*}_{e}(A,q,\eta)$ and thus Theorem $\ref{thmexistence}$ holds for all $c\geq c^{*}_{e}(g)$. This is one way to understand where does the threshold $c^{*}_{e}(A,q,\eta)$ comes from.

\bigskip

For positive reaction terms, one has $\lambda_{1}'\leq 0$ but not necessarily $\lambda_{1}'<0$. It is possible to prove that the previous result is still true when $\lambda_{1}'=0$, which is new even if only the advection is heterogeneous. 
\begin{thm} \label{thmexistencelambda0}
Assume that $f(t,x,0)=f(t,x,1)=0$, that $f(t,x,s)>0$ for all $s\in (0,1)$ and that $f(t,x,s)<0$ if $s>1$. Then for all unit vector $e$, for all speed $c\geq c^{*}_{e}(A,q,\eta)$, where $\eta(t,x)=\sup_{0<s<p(t,x)}\frac{f(t,x,s)}{s}$, there exists a pulsating traveling front $u$ of speed $c$ in direction $-e$ that links $0$ to $1$. The profile $\phi$ is nondecreasing almost everywhere in $z$.

If $c<c^{*}_{e}(A,q,\mu)$, where $\mu(t,x)=f_u'(t,x,0)$, then there exists no pulsating traveling fronts of speed $c$.
\end{thm}

This kind of existence result for flat nonlinearity has been proved before in the case of cylinders with orthogonal heterogeneity by Berestycki and Nirenberg \cite{Nirenberg} and in space periodic environments by Berestycki and Hamel \cite{Frontexcitable}. 

This proves that our hypothesis $\lambda_{1}'<0$ is not optimal for some particular $f$. But there exists a threshold for the non-existence of pulsating traveling front for general $f$:
\begin{prop} \label{cexlambda1'=0}
If $k_{0}(A,q,\eta)>0$, then there exists no bounded positive entire solution of equation $(\ref{eqprinc})$. In particular, there exists no pulsating traveling front.
\end{prop}
It is not clear what happens between the thresholds $k_{0}(A,q,\mu)<0$ and $k_{0}(A,q,\eta)>0$.

\bigskip

We can now define the $minimal$ $speed$ of the pulsating traveling fronts for general nonlinearities:
\begin{thm} \label{mainresult}
Assume that $\lambda_{1}'< 0$ and that Hypothesis $\ref{uniquenesshyp}$ is satisfied. Then for all unit vector $e$, there exists a minimal speed $c^{*}_{e}$ such that there exists a pulsating traveling front $u$ of speed $c^{*}_{e}$ in direction $-e$ that links $0$ to $p$, while no such front exists if $c<c^{*}_{e}$. 
This minimal speed satisfies: 
\begin{equation} \label{estimatec} c^{*}_{e}(A,q,\mu)\leq c^{*}_{e}\leq c^{*}_{e}(A,q,\eta), \end{equation}
where $\eta(t,x)=\sup_{0<s<p(t,x)}\frac{f(t,x,s)}{s}$ and $\mu(t,x)=f_u'(t,x,0)$. 
\end{thm}

The classical result for homogeneous environments is that $for$ $all$ $c\geq c^{*}_{e}$, there exists a pulsating traveling front. This result has been extended to space periodic environment in \cite{Base2}. If $f$ is of KPP type, then Theorem $\ref{thmexistence}$ yields that this is true since $c^{*}_{e}=c^{*}_{e}(A,q,\mu)$. This is not clear if this result still holds if $f$ is not of KPP type

If $f$ is not of KPP type, set 
\[\mathcal{C}=\{c\in\R, \ \hbox{there exists some pulsating traveling front of speed} \ c\}.\]
We do not know if $\mathcal{C}$ is a half-line, like in the homogeneous case, but we still know from Theorem $\ref{thmexistence}$ that there exists at least a half line $[c^{*}_{e}(A,q,\eta),+\infty)$ included in $\mathcal{C}$ and a half-line $(-\infty,c^{*}_{e}(A,q,\eta))$ in its complementary. Moreover, Theorem $\ref{mainresult}$ is in fact a corollary of the following proposition:
\begin{prop} \label{vitferme}
Assume that Hypothesis $\ref{uniquenesshyp}$ is satisfied and that $\mathcal{C}$ is not empty. Then the set $\mathcal{C}$ is closed.
\end{prop}

\subsection{Regularity of the fronts}

The next theorem is the most complete theorem of this paper and it answers to all the open questions we pointed out in the introduction. Namely, it states that one can get $Lipschitz \ continuous$ pulsating traveling fronts (and not only measurable ones) $for \ all$ $c\geq c^{*}_{e}$. But one particular hypothesis is needed, which is somehow a strong global version of the classical KPP hypothesis: $s\mapsto f(t,x,s)/s$ is nonincreasing. The previous theorems gave partial answers to the open questions, but under weaker hypotheses.

\begin{thm} \label{mainresultcont}
Assume that $\lambda_{1}'< 0$, that Hypothesis $\ref{uniquenesshyp}$ is satisfied and that $s\mapsto f(t,x,s)/s$ is nonincreasing for all $(t,x)\in\R\times\R^{N}$. Then for all $c\geq c^{*}_{e}$, there exists a $Lipschitz \ continuous$ pulsating traveling front of speed $c$, while there exists no such front for $c<c^{*}_{e}$.
\end{thm}

\noindent {\bf Remark:} The reader may remark that the hypotheses of this theorem are satisfied in particular if $\lambda_{1}'<0$ and $(\ref{hyp1})$ and $(\ref{hyp2})$ are satisfied. Furthermore, as $s\mapsto f(t,x,s)/s$ is nonincreasing, then $f$ is of KPP type and $c^{*}_{e}=c^{*}_{e}(A,q,\mu)$.

\bigskip

If $f$ does not satisfy the hypotheses of Theorem $\ref{mainresultcont}$, the method that is used in the proof of this theorem fails, even when $c>c^*_e(A,q,\eta)$. There is no particular heuristic reason why this theorem might be false if these hypotheses are not satisfied and it only seems to be a technical issue. Anyway, we can notice that a pulsating traveling front $v=v(y,t,x)$ in the sense of Definition $\ref{deftravelingfronts}$ only admits a countable number of points of discontinuity in $y$. Namely, one can assume, up to some modifications on a set of measure $0$, that $v$ is nondecreasing in $z$. Thus $v$ admits a limit on its left and on its right everywhere, which yields that the discontinuity points are isolated and thus their set is countable.

\subsection{Spreading properties}

It is not possible to prove that there exists no pulsating traveling fronts for small speeds using the classical methods. Thus, we had to prove spreading properties for front-like initial data in order to get the non-existence result. We used the same method as Mallordy and Roquejoffre \cite{Roquejoffre} and Nolen, Rudd and Xin \cite{Nolen}. 

The following results have already been proved by Weinberger \cite{Weinberger} in a time and space discrete context. Considering the Poincar\'e map of time $T$ that is associated with equation $(\ref{eqprinc})$, the following results are then consequences of the results of Weinberger. In \cite{posspreadspeed}, we used this method in order to investigate the case of compactly supported initial data and to give an alternative and independent proof to that of Weinberger. Furthermore, we managed to use this method in a general heterogeneous media and to prove the existence of a positive spreading speed when $f$ is positive and $q\equiv 0$. 

\begin{prop} \label{spreadingsub}
Assume that $\lambda_{1}'\leq 0$ and that Hypothesis $\ref{uniquenesshyp}$ is satisfied. Take $u_{0}\leq p$ a nonnegative continuous initial datum and an interval $[a_{1},a_{2}]\subset \R$ such that:
\[\inf_{x\in\R^{N}, e\cdot x \in [a_{1},a_{2}]}u_{0}(x)>0 .\]
Then for all $c<c^{*}_{e}(A,q,\mu)$, the solution $u$ of equation $(\ref{eqprinc})$ associated with the initial datum $u_{0}$ satisfies:
\[u(t,x-cte)-p(t,x-cte)\rightarrow 0 \ \hbox{as} \ t\rightarrow +\infty, \]
locally uniformly in $x\in\R^{N}$.
\end{prop}

If the speed $c$ is larger than the speed $c^{*}_{e}(A,q,\eta)$, we get the opposite spreading property:

\begin{prop} \label{spreadingsup}
Assume that $\lambda_{1}'\leq 0$, that Hypothesis $\ref{uniquenesshyp}$ is satisfied and that the growth rate $\eta:s\mapsto f(t,x,s)/s$ is bounded from above for all $(t,x)\in\R\times\R^{N}$. Take $u_{0}$ a nonnegative continuous bounded initial datum and assume that $u_{0}(x)=O_{x\cdot e \rightarrow -\infty}(e^{\lambda x\cdot e})$ for all $0<\lambda<\lambda_{c^{*}(\eta)}(\eta)$. 
Under these hypotheses, for all $c>c^{*}_{e}(A,q,\eta)$, the solution $u$ of equation $(\ref{eqprinc})$ associated with the initial datum $u_{0}$ satisfies:
\[u(t,x-cte)\rightarrow 0 \ \hbox{as} \ t\rightarrow +\infty,\]
uniformly in $x\in \{x\in\R^{N}, \ x\cdot e\leq B\}$ for all $B\in\R$.
\end{prop}

\noindent {\bf Remark:} The quantity $\lambda_{c^{*}(\eta)}(\eta)$ will be defined in section $\ref{deflambdac}$, but one may already remark that any initial datum with compact support satisfies the hypotheses of this proposition.


\section{Preliminaries}

This section is devoted to the proof of the equivalence between Definitions $\ref{deftravelingfrontsNolen}$ and $\ref{deftravelingfronts}$ and to some technical results on the family of periodic principal eigenvalues $(k_\lambda)_\lambda$. 

\subsection{Equivalence of the definitions}

We first prove that the two Definitions $\ref{deftravelingfrontsNolen}$ and $\ref{deftravelingfronts}$ are equivalent.

Assume first that $\phi$ satisfies the properties of Definition $\ref{deftravelingfrontsNolen}$ and set:
\[v(y,t,x)=\phi(y+x\cdot e+ct,t,x).\]
This function satisfies:
\begin{equation} \label{ueps}
 \partial_{t}v-\nabla\cdot (A\nabla v)+q\cdot\nabla v=f(t,x,v) \ \hbox{in} \ \mathcal{D}'(\R\times\R\times\R^{N}).
\end{equation}

Take $\theta\in W^{2,2}(\R\times\R^{N})$ a compactly supported function that only depends on $t$ and $x$. Set:
\begin{equation} \begin{array}{rcl}
 g_{\theta}(y)&=&<\partial_{t}v-\nabla\cdot (A\nabla v)+q\cdot\nabla v=f(t,x,v),\theta>_{\mathcal{D}'(\R\times\R^{N})\times\mathcal{D}(\R\times\R^{N})}\\
 &&\\
&=&\int_{\R\times\R^{N}}\big(v(\partial_{t}\theta-  \nabla\cdot (A\nabla \theta)+ \nabla\cdot (q \theta))-f(t,x,v)\theta \big) dtdx.\\
\end{array}
\end{equation}
This function belongs to $L^{1}_{loc}(\R)$ since $v\in L^{\infty}(\R\times\R\times\R^{N})$. Moreover, equation $(\ref{ueps})$ yields that for all compactly supported function $\chi\in L^{\infty}(\R)$, one has:
\[\int_{\R}g_{\theta}(y)\chi(y)dy=0.\]
Thus, there exists some set $\mathcal{N}_{\theta}$ of null measure such that for all $y\notin\mathcal{N}_{\theta}$, one has $g_{\theta}(y)=0$. 

In the other hand, set:
\[\mathcal{H}=\{ h\in W^{2,2}(\R\times\R^{N}), h \ \hbox{is compactly supported}\}.\]
One knows that there exists a sequence $(\theta_{n})_{n}\in\mathcal{H}$ which is dense in $\mathcal{H}$. The set $\mathcal{N}=\cup_{n\in\N}\mathcal{N}_{\theta_{n}}$ is of null measure. 

As for all $y\notin\mathcal{N}$ and for all $n$, one has $g_{\theta_{n}}(y)=0$, one gets for all $R>0$:
\[\int_{B_{R}\cap \mathcal{N}^{c}}|g_{\theta_{n}}|dy = 0.\]
Next, one can remark that the linear function $\theta \in \mathcal{H}\mapsto g_{\theta}\in L^{1}_{loc}(\R)$ is continuous and thus for all $R>0$ and for all $\theta\in \mathcal{H}$:
\[\int_{B_{R}\cap \mathcal{N}^{c}}|g_{\theta}|dy = 0,\]
which means that there exists some null measure set $\mathcal{M}$ such that for all $y\notin\mathcal{M}\cup\mathcal{N}$, for all $\theta\in \mathcal{H}$, one has $g_{\theta}(y)=0$, which means that $g_{\theta}=0$ almost everywhere.

Thus for almost every $y\in\R$, the function $v(y,\cdot,\cdot)$ is a weak solution of:
\begin{equation}
 \partial_{t}v-\nabla\cdot (A\nabla v)+q\cdot\nabla v=f(t,x,v) \ \hbox{in} \ \R\times\R^{N}.
\end{equation}
The Schauder parabolic estimates give that $u: (t,x)\mapsto v(0,t,x)$ satisfies the hypotheses of Definition $\ref{deftravelingfronts}$. 

\bigskip

In the other hand, assume that $u$ satisfies the hypotheses of Definition $\ref{deftravelingfronts}$ and set $v(y,t,x)=\phi (y+x\cdot e+ct,t,x)$. Take any family of compactly supported functions $(\chi_{k})_{k\in [1,n]}\in L^{\infty}(\R)$ and $(\theta_{k})_{k\in [1,n]}\in \mathcal{D}(\R\times\R^{N})$. One has for almost every $y\in\R$:
\[\chi_{k}(y)\theta_{k}(t,x) \big( \partial_{t}v-\nabla\cdot (A\nabla v)+q\cdot \nabla v -f(t,x,v))=0.\]
Thus, integrating over $\R\times\R\times\R^{N}$, one gets:
\[\begin{array}{l}<\partial_{t}v-\nabla\cdot (A\nabla v)+q\cdot \nabla v -f(t,x,v),\sum_{k=1}^{n}\chi_{k}\theta_{k}>_{\mathcal{D}'\times\mathcal{D}}\\
=\sum_{k=1}^{n}\int_{R\times\R\times\R^{N}}v(y,t,x) \big( \partial_{t}\chi_{k}\theta_{k}-\nabla\cdot (A\nabla \chi_{k}\theta_{k})+\nabla\cdot(q\theta_{k}\chi_{k})\big)\\
-\sum_{k=1}^{n}\int_{R\times\R\times\R^{N}}\chi_{k}\theta_{k} f(t,x,v)\\
=0.\\ \end{array}\]

As the regular functions with separated variables span a dense subset of the compactly supported $L^{\infty}(\R, W^{2,2}(\R\times\R^{N}))$ functions, one finally gets 
\[\partial_{t}v-\nabla\cdot (A\nabla v)+q\cdot \nabla v -f(t,x,v) \ \hbox{in} \ \mathcal{D}'(\R\times\R\times\R^{N}).\]
Thus the function $\phi(z,t,x)=v(z-x\cdot e-ct,t,x)$ satisfies the hypotheses of Definition $(\ref{deftravelingfrontsNolen})$.$\Box$

\subsection{More properties of the eigenvalue family}

In order to build an invariant domain, we first need to prove or recall some more precise properties for the family of eigenvalues $(k_{\lambda e})_{\lambda>0}$.

\begin{prop} \cite{eigenvalue}\label{concavitykalpha}
Set $F$ the map:
\[ \begin{array}{ccl}
\mathbb{R}^{N} \times \mathcal{C}_{per}^{0}(\mathbb{R}^{N} \times \mathbb{R}) &\rightarrow& \mathbb{R}\\
(\lambda,\mu) &\mapsto& k_{\lambda}(\mu)\\
\end{array} \]

Then:
\[1) F(\lambda,\mu)= \max_{\phi \ is \ periodic \ in \ (t,x), \ \phi>0.}\min_{\mathbb{R} \times \mathbb{R}^{N}} \Big( \frac{L_{\lambda} \phi }{\phi} \Big)=\min_{\phi \ is \ periodic \ in \ (t,x), \ \phi>0.}\max_{\mathbb{R} \times \mathbb{R}^{N}} \Big( \frac{L_{\lambda} \phi }{\phi} \Big).\]

2) $F$ is concave and continuous.

3) There exists $\beta \in \mathbb{R}$ such that for all $\lambda \in\mathbb{R}^{N}$: \begin{equation} \label{ineq}-\|\mu\|_{\infty} - \beta \|\lambda\|- \Gamma\|\lambda\|^{2}\leq k_{\lambda}(\mu) \leq \|\mu\|_{\infty} + \beta \|\lambda\|- \gamma\|\lambda\|^{2}\end{equation}
where $\gamma$ and $\Gamma$ have been defined in $(\ref{ellipticity})$.

4) For all $\mu$, $F(.,\mu)$ reaches its maximum on $\mathbb{R}^{N}$. 
\end{prop}

\smallskip

These properties enable us to define and characterize a quantity
that we will use to get good bounds for the propagation speed.
Namely, for all $(A,q,\mu)$, $e\in\mathbb{S}^{N-1}$ and for all $\gamma>\epsilon\geq 0$, we set:
\[ c^{*}_{e,\epsilon}(A,q,\mu)=\min \{ c\in \R, \ \hbox{there exists} \ \lambda>0 \
\hbox{such that} \ k_{\lambda e}(A,q,\mu)+\lambda c-\epsilon\lambda^{2}=0 \}.\] In the
sequel, $A,q$ and $e$ will be fixed and we will denote
$c^{*}_{\epsilon}(\mu)=c^{*}_{e,\epsilon}(A,q,\mu)$.

\begin{prop} \label{deflambdac} Take $\gamma>\epsilon>0$. 

1) For $c=c^{*}_{\epsilon}(\mu)$, there exists exactly one solution $\lambda_{c}^{\epsilon}(\mu)>0$ to equation \[k_{\lambda e}(A,q,\mu)+\lambda c-\epsilon\lambda^{2}=0.\]

2) For all $c>c^{*}_{\epsilon}(\mu)$, there exist exactly two $\lambda>0$ that
satisfy $k_{\lambda e}(A,q,\mu)+\lambda c-\epsilon\lambda^{2}=0$. We denote
$\lambda_{c}^{\epsilon}(\mu)\leq \Lambda_{c}^{\epsilon}(\mu)$ those two solutions
$\lambda$.

3) One has the characterization:
\[c^{*}_{\epsilon}(\mu)=\min_{\lambda>0}\frac{-k_{\lambda e}(A,q,\mu)+\epsilon\lambda^{2}}{\lambda}.\]
\end{prop}

\noindent {\bf Proof.}
1) Set $G(c)=\max_{\lambda>0}(k_{\lambda e}(A,q,\mu)+\lambda c-\epsilon\lambda^{2})$. This maximum is reached because of estimate 2) of Proposition $\ref{concavitykalpha}$ and thus this function is increasing in $c$. If $c<\min_{\lambda>0}\frac{-k_{\lambda e}(\mu)+\epsilon\lambda^{2}}{\lambda}$, one easily gets $G(c)<0$. In the other hand, (\ref{ineq}) yields that:
\[G(c)\geq -\|\mu\|_{\infty}+\max_{\lambda>0}((c-\beta)\lambda-\lambda^{2}(\Gamma+\epsilon))=-\|\mu\|_{\infty}+\frac{(c-\beta)^{2}}{4(\Gamma+\epsilon)}\]
if $c\geq \beta$, which gives $G(c)\rightarrow+\infty$ as $c\rightarrow+\infty$. Thus, $c=c^{*}_{\epsilon}(\mu)$ is the only solution to equation $F(c)=0$ 

Next, fix $c=c^{*}_{\epsilon}(\mu)$. The Kato-Rellich perturbation theorem (see \cite{Kato}) yields that the function $\lambda\mapsto k_{\lambda e}(A,q,\mu)$ is analytic with respect to $\lambda$, locally uniformly in $(A,q,\mu)$. Thus, if there exist two solutions $\lambda_{1}<\lambda_{2}$ of equation $k_{\lambda e}(A,q,\mu)+\lambda c -\epsilon\lambda^{2}= 0$, then we know from the definition of $c^*_\epsilon(\mu)$ and the concavity of $\lambda\mapsto k_{\lambda e}(A,q,\mu)$ that $k_{\lambda e}(A,q,\mu)+\lambda c-\epsilon\lambda^{2} = 0$ for all $\lambda\in [\lambda_{1},\lambda_{2}]$. The isolated zeros principle would give $k_{\lambda e}(A,q,\mu)+\lambda c-\epsilon\lambda^{2} = 0$ for all $\lambda\in\R$ in this case. As $ k_{0}(A,q,\mu)=\lambda_{1}'(A,q,\mu)<0$, this gives a contradiction. Thus for $c=c_{\epsilon}^{*}(\mu)$, there exists a unique $\lambda_{c}^{\epsilon}(\mu)$ such that $k_{\lambda_{c}(\mu) e}(A,q,\mu)+\lambda_{c}(\mu) c -\epsilon\lambda^{2}= 0$.

2) If $c>c^{*}_{\epsilon}(\mu)$, the previous step yields that the maximum of the function $\lambda\mapsto k_{\lambda e}(A,q,\mu)+\lambda c-\epsilon\lambda^{2}$ is positive. As this function is concave, negative for $\lambda=0$ and goes to $-\infty$ as $\lambda\rightarrow+\infty$, one easily concludes.

3) This easily follows from the characterization $G(c_{\epsilon}^{*}(\mu))=0$. $\Box$

\bigskip

Lastly, we need a continuity property for $\epsilon\mapsto \lambda_{\epsilon}^{c}$:

\begin{prop} \label{contlambdaeps}
If $\epsilon_{n}\rightarrow \epsilon\geq 0$, where $c>c_{\epsilon_{n}}$ for all $n$, then one has:
\[\lambda_{\epsilon_{n}}^{c}\rightarrow \lambda_{\epsilon}^{c} \ \hbox{and} \ \Lambda_{\epsilon_{n}}^{c}\rightarrow \Lambda_{\epsilon}^{c}\]
as $n\rightarrow+\infty$.
\end{prop}

\noindent {\bf Proof.}
Using estimate $1)$ of Proposition $\ref{concavitykalpha}$, one knows that the sequences $(\lambda_{\epsilon_{n}}^{c})_{n}$ and $(\Lambda_{\epsilon_{n}}^{c})$ are bounded. One can extract two converging subsequences $\lambda_{\epsilon_{n'}}^{c}\rightarrow \lambda_{\infty}^{c}$ and $\Lambda_{\epsilon_{n'}}^{c}\rightarrow \Lambda_{\infty}^{c}$. The continuity of $\lambda\mapsto k_{\lambda e}(A,q,\mu)$ yields that:
\[k_{\lambda_{\infty}^{c} e}(A,q,\mu)+\lambda_{\infty}^{c} c-\epsilon(\lambda_{\infty}^{c})^{2}=0, k_{\Lambda_{\infty}^{c} e}(A,q,\mu)+\Lambda_{\infty}^{c} c-\epsilon(\Lambda_{\infty}^{c})^{2}=0 \ \hbox{and} \ \lambda_{\infty}^{c}\leq \Lambda_{\infty}^{c}.\]
The previous proposition thus gives $\lambda_{\infty}^{c}=\lambda_{\epsilon}^{c}$ and $\Lambda_{\infty}^{c}=\Lambda_{\epsilon}^{c}$. $\Box$


\section{The KPP case}

This section is devoted to the proof of Theorem $\ref{thmexistence}$ and characterization $(\ref{characterizationminspeed})$. We will directly prove that there exists a pulsating traveling front of speed $c$ for all $c> c^{*}_{e}(A,q,\mu)$ and a result that will be proved later will give the existence of a pulsating traveling front of speed $c^{*}_{e}(A,q,\mu)$. We use the same kind of invariant domains as in \cite{BHKR, GuoHamel1, HamelRyzhik}. 

\subsection{Study of the regularized problem in finite cylinders}

In all this section, we fix $c> c^{*}(\mu)$. Equation $(\ref{eqtf})$ exhibits two main issues. First of all, it is defined in an infinite domain in $z$. Secondly, it is a degenerate parabolic equation. We will first solve a modified regular parabolic equation in a finite domain. Then, we will pass to the limit. 

Let $\Sigma_{a}$ be the domain $(-a,a)\times\mathbb{R}\times \mathbb{R}^{N}$. We first investigate the equation:
\begin{equation} \label{regfinite0}
\left\{ \begin{array}{l}
L_{\epsilon} \phi = f(t,x,\phi) \hbox{ in } \Sigma_a, \\
\phi \ \hbox{is periodic in} \ t,x \\
\end{array}
\right.
\end{equation}
where $L_{\epsilon}$ is the regular parabolic operator defined by:
\begin{equation} \begin{array}{rcl}
L_{\epsilon} \phi &=& \partial_{t} \phi - \nabla \cdot (A(t,x) \nabla \phi) - (eA(t,x)e+\epsilon) \partial_{zz}\phi - \nabla \cdot (A(t,x)e \partial_{z} \phi)\\
&&-\partial_{z}(e A(t,x)\nabla \phi) 
+q(t,x)\cdot\nabla\phi+q(t,x)\cdot e\partial_{z}\phi+c\partial_{z}\phi.\\
\end{array}
\end{equation}
The condition in $z=\pm a$ will be fixed later.

We first construct a subsolution which does not depend on $a$. This kind of subsolution has first been used in \cite{BHKR}. This gives some lower bound for the solution of $(\ref{regfinite0})$ that enables us to easily pass to the limit $a\rightarrow +\infty$. This lemma is one of the main point of our proof: as it is available for all $c\geq c^{*}(\mu)$, it directly proves in the KPP case that the speed we will finally obtain is minimal.

In \cite{Nolen}, the authors used a different approach. They first proved the existence of pulsating traveling fronts for ignition type nonlinearities, that is, there exists some $\theta\in (0,1)$ such that $f(s)=0$ if $s\in (0,\theta)$, $f(s)>0$ if $s\in (\theta,1)$ and $f(1)=0$. In this case, there exists a unique speed associated with some pulsating traveling front. Then, they let $\theta$ go to zero and get some pulsating traveling front for $one$ speed $c^{*}$, which is the limit of the previous speeds. The hard part is then to prove that there exists a pulsating traveling front for all $c>c^{*}$. In the present paper, we overcome this issue by directly proving this existence result for a half-line of speeds.

This subsolution exists if $f$ satisfies the following hypothesis: there exists some $\beta>0$, $\rho>0$ and $r>0$ such that
\begin{equation} \label{thetafC2} \forall 0\leq u < \beta, \ \forall t,x, \ \mu(t,x)u \leq \rho u^{1+r} +f(t,x,u). \end{equation}
This is true in particular if $f\in\mathcal{C}^{1,r}(\R\times\R^{N}\times [0,\beta])$, which is one of our hypotheses, but the reader can check that the results of this article hold if only $(\ref{thetafC2})$ is satisfied. This condition is a little bit sharper than the $\mathcal{C}^{1}$ regularity for $f$. A linear coupled equation was considered in \cite{BHKR} and this kind of condition has first been stated for a nonlinear equation in \cite{GuoHamel1}. It has also been used by Bag\`es \cite{TheseBages} to prove that, in space periodic media, one can construct some pulsating traveling front with a given exponential decay at infinity. 

Bag\`es also managed to construct a subsolution in the critical case $c=c^*_e$. The kind of subsolution he constructed cannot be used here since we consider a regularized problem. The subsolution available for $c=c^*_e$ is not a subsolution of the regularized problem anymore. Thus we will first consider the case $c>c^*_e$ and then we will pass to the limit $c\rightarrow c^*_e$. 

\begin{lem} \label{thetalem}
For all $c> c^{*}_{\epsilon}(\mu)$, set 
\begin{equation} \theta_{0,\epsilon}(z,t,x)=\psi_{\mu,\lambda^{\epsilon}_{c}}(t,x)e^{\lambda^{\epsilon}_{c}z}-A\psi_{\mu,\lambda^{\epsilon}_{c}+\gamma}
(t,x) e^{(\lambda^{\epsilon}_{c}+\gamma)z}
\end{equation} 
where the functions $\psi_{\mu,\lambda}$ are eigenfunctions associated with $L_{\lambda}$ and normalized by $\|\psi_{\mu,\lambda}\|_{\infty}=1$. Then there exist $\gamma, A$ which do not depend on $\epsilon$ such that  the function $\theta_{\epsilon}=\max \{\theta_{0,\epsilon}, 0\}$ satisfies $\theta_{\epsilon}\leq p$ and: 
\begin{equation}
L_{\epsilon}\theta_{\epsilon}\leq 0 \ \hbox{in the sense of distribution}.
\end{equation}
\end{lem}

\noindent {\bf Proof.}
Set $\alpha_{\lambda}^{c,\epsilon}=k_{\lambda e}(\mu)+\lambda c-\epsilon\lambda^{2}$. As $ c> c^{*}_{\epsilon}(\mu)$, one knows from Proposition \ref{deflambdac} that the equation $\alpha^{c,\epsilon}_{\lambda}=0$ only admits one or two solutions $\lambda_{c}^{\epsilon} < \Lambda_{c}^{\epsilon}$ for all $\epsilon\geq 0$. As $\epsilon\geq 0\mapsto \Lambda_{c}^{\epsilon} - \lambda_{c}^{\epsilon}$ is continuous, it admits a positive infimum $h>0$. As $\lambda \mapsto \alpha_{\lambda}^{c,\epsilon}$ is concave, one gets that $\alpha^{c,\epsilon}_{\lambda_{c}^\epsilon+\gamma}>0$ if $\gamma<h$. We also choose $\gamma$ independent of $\epsilon$ such that $\lambda_{c}^{\epsilon}+\gamma\leq (1+r)\lambda_{c}^{\epsilon}$. 

We know from $(\ref{thetafC2})$ that:
\[ \forall \ 0\leq u < \beta, \ \forall t,x, \ \mu(t,x)u \leq \rho u^{1+r} +f(t,x,u). \]

We can fix $A$ independent of $\epsilon$ such that:
\[\max_{(z,t,x)\in\R\times\R\times\R^{N}}\theta_{0,\epsilon}(z,t,x)\leq \min \{\beta,p\},\]
\[\forall (t,x)\in\R\times\R^{N}, \rho \psi_{\lambda_{c}^{\epsilon}}^{1+r}(t,x)\leq A\alpha^{c,\epsilon}_{\lambda_{c}^{\epsilon}+\gamma}\psi_{\lambda_{c}^{\epsilon}+\gamma},\]
and such that for all $(z,t,x)\in\R^{+}\times\R\times\R^{N}$, one has $\theta_{0,\epsilon}(z,t,x)\leq 0$.

Now that it is proved that $A$ and $\gamma$ can be chosen independent of $\epsilon$, we forget the dependence in $\epsilon$ in the notations in order to make the proof easier to read. 
Set: $\Omega^{+}=\{(z,t,x) \in \mathbb{R}\times\mathbb{R}\times\mathbb{R}^{N} ; \theta_{0}(z,t,x) > 0 \}$. We compute for all $(z,t,x) \in \Omega^{+}$:
\[
\begin{array}{rcl}
 L_{\epsilon} (\theta_{0})&=&\mu(t,x) \theta_{0}-A\alpha^{c,\epsilon}_{\lambda_{c}^{\epsilon}+\gamma}\psi_{\lambda_{c}^{\epsilon}+\gamma}(t,x) e^{(\lambda_{c}^{\epsilon}+\gamma)z} \\
&& \\
 &\leq& f(t,x,\theta_{0})+\rho \theta_{0}^{1+r}-A\alpha^{c,\epsilon}_{\lambda_{c}^{\epsilon}+\gamma}\psi_{\lambda_{c}^{\epsilon}+\gamma} e^{(\lambda_{c}^{\epsilon}+\gamma)z}\\
&&\\
 &\leq&
 f(t,x,\theta_{0})+\rho \psi_{\lambda_{c}^{\epsilon}}^{1+r}e^{\lambda_{c}^{\epsilon}(1+r)z}-A\alpha^{c,\epsilon}_{\lambda_{c}^{\epsilon}+\gamma}\psi_{\lambda_{c}^{\epsilon}+\gamma} e^{(\lambda_{c}^{\epsilon}+\gamma)z}\\
&&\\
 &\leq&
 f(t,x,\theta_{0})+\big(\rho \psi_{\lambda_{c}^{\epsilon}}^{1+r}-A\alpha^{c,\epsilon}_{\lambda_{c}^{\epsilon}+\gamma}\psi_{\lambda_{c}^\epsilon+\gamma}\big) e^{(\lambda_{c}^{\epsilon}+\gamma)z}\\
&&\\
 &\leq&
 f(t,x,\theta_{0}).\\
\end{array} \]
Thus, $\theta_0$ is a subsolution of equation $(\ref{regfinite0})$ over $\Omega^{+}$. As $0$ is a solution of equation $(\ref{regfinite0})$, the Hopf lemma gives that $\theta=\max \{\theta_{0},0\}$ is a subsolution of equation $(\ref{regfinite0})$ in the sense of distributions. $\Box$

\bigskip

\begin{lem}
For all $\epsilon\in (0,\epsilon_{0})$, the function 
\[\zeta(z,t,x)=\min \{ 0, \psi_{\mu,\lambda_{c}^{\epsilon}(\mu)}(t,x)e^{\lambda_{c}^{\epsilon}(\mu)z}\}\]
is a supersolution of equation $(\ref{regfinite0})$ in the sense of distributions.
\end{lem}

\noindent {\bf Proof.}
One easily computes:
\begin{equation} \begin{array}{rcl}
L_{\epsilon}(\psi_{\lambda_{c}^{\epsilon}} e^{\lambda_{c}^{\epsilon}z}) &=& \big\{ \partial_{t}  \psi_{\lambda_{c}^{\epsilon}} - \nabla\cdot(A\nabla \psi_{\lambda_{c}^{\epsilon}}) - 2\lambda_{c}^{\epsilon} e A\nabla \psi_{\lambda_{c}^{\epsilon}}+q\cdot\nabla \psi_{\lambda_{c}^{\epsilon}}\\
&&- (\lambda_{c}^{\epsilon}\nabla\cdot(A e)+(1+\epsilon)(\lambda_{c}^{\epsilon})^{2}+\lambda_{c}^{\epsilon} c+\lambda_{c}^{\epsilon}q\cdot e )\psi_{\lambda_{c}^{\epsilon}}^{\epsilon} \big\} e^{\lambda_{c}^{\epsilon}z}\\
&=&\mu(t,x)\psi_{\lambda_{c}^{\epsilon}} e^{\lambda_{c}^{\epsilon}z}
\geq f(t,x,\psi_{\lambda_{c}^{\epsilon}} e^{\lambda_{c}^{\epsilon}z})\\
\end{array}\end{equation}

As the minimum of two supersolutions is a supersolution in the sense of distributions, this gives the conclusion. $\Box$

\bigskip

We are now able to define the boundary conditions associated with equation $(\ref{regfinite0})$. Namely, we solve:
\begin{equation} \label{regfinite}
\left\{ \begin{array}{l}
L_{\epsilon} \phi = f(t,x,\phi) \ \hbox{in} \ \Sigma_{a},\\
\phi \ \hbox{is periodic in} \ t,x, \\
\phi (-a,t,x)= \theta_{\epsilon} (-a,t,x) \ \hbox{for all} \ (t,x)\in\R\times\R^{N},\\
\phi (a,t,x)= \zeta_{\epsilon}(t,x) \ \hbox{for all} \ (t,x)\in\R\times\R^{N}.\\
\end{array}
\right.
\end{equation}

The boundary condition used in $-a$ may seem odd and one can wonder why we use this complicated condition instead of $0$ for example, like in \cite{Base2}. In fact, this condition will enable us to put the subsolution $\theta_{\epsilon}$ under the solution of equation $(\ref{regfinite})$.

\begin{lem} \label{existencerefinite}
For all $\epsilon>0$, there exists a strong solution $\phi_{a}^{\epsilon}$ of $(\ref{regfinite})$ in $\mathcal{C}^{2,1,2}(\Sigma_{a})$ that satisfies
\[\forall (z,t,x)\in\Sigma_{a}, \ \theta_{\epsilon}(z,t,x)\leq \phi_{a}^{\epsilon}(z,t,x)\leq \zeta_{\epsilon}(z,t,x).\]
Moreover, there exists some $a_{0}$ such that for all $a>a_{0}$, the function $\phi_{a}^{\epsilon}$ is nondecreasing with respect to $z$.
\end{lem}

\noindent {\bf Proof.}
We know that $\theta_{\epsilon}$ and $\zeta_{\epsilon}$ are sub and supersolutions of equation $(\ref{regfinite})$. 
Let $r$ be the function defined by:
\[r(z,t,x)=\zeta_{\epsilon}(t,x)\frac{z+a}{2a}+\theta_{\epsilon}(z,t,x)\frac{z-a}{2a}.\]
We make the change of variable $u=\phi-r$. Then $\phi$ satisfies $(\ref{regfinite})$ if and only if $u$ satisfies:
\begin{equation} \label{finite} \left\{ \begin{array}{l}
L_{\epsilon}u=g(t,x,u) \ \hbox{in} \ \Sigma_{a},\\
u \ \hbox{is} \ T,L_{1},...,L_{N}- \hbox{periodic},\\
u(-a,.,.)=0, \ u(a,.,.)=0, \\
\end{array} \right. \end{equation}
where:
\begin{equation} \begin{array}{l}g(t,x,u)=f(t,x,u+r)-L_{\epsilon}r.\end{array}\end{equation}
This function is locally Lipschitz continuous with respect to $u$.

As $\theta_{\epsilon}-r$ is a subsolution and $\zeta_{\epsilon}-r$ is a supersolution of equation $(\ref{finite})$, in order to conclude using an iteration procedure, we need to prove that for $\beta$ sufficiently large, the operator $L_{\epsilon}+\beta$ is invertible. Take some $\beta>\frac{1}{2}\|\nabla\cdot q\|_{\infty}$ and set $\mu=\beta-\frac{1}{2}\|\nabla\cdot q\|_{\infty}$.

Take some $g\in\mathcal{C}^{0}(\Sigma_{a})$ such that $g$ is periodic in $t$ and $x$ and
\[g(-a,t,x)=g(a,t,x)=0 \ \hbox{for all} \ (t,x)\in\R\times\R^{N}.\]
Set $\mathcal{L}=\{u_{0}\in L^{2}((-a,a)\times C, u \ \hbox{is L-periodic}, \ u(-a,x)=u(a,x)=0)\}$ (this set is in fact the closure of the set of the continuous functions that satisfy the boundary conditions with respect to the $L^{2}$ norm).
For all $u_{0}\in \mathcal{L}$, we define $(z,t,x)\mapsto u(z,t,x)\in\mathcal{C}^{1}(\R^{+},\mathcal{L})$ the solution of 
\begin{equation} \left\{ \begin{array}{l}
L_{\epsilon}u+\beta u=g,\\
u(z,0,x)=u_{0}(z,x),\\
\end{array} \right. \end{equation}
and we investigate the map 
\begin{equation} \left\{ \begin{array}{rcl}
G: \mathcal{L} &\rightarrow & \mathcal{L}\\
u_{0}& \mapsto & u(T).\\
\end{array} \right. \end{equation}

Take $u_{1},u_{2}\in\mathcal{L}$ and set $U(z,t,x)=(u_{1}(z,t,x)-u_{2}(z,t,x))e^{\mu t}$. This function satisfies:
\[\begin{array}{l}\partial_{t} U - \nabla \cdot (A(t,x) \nabla U) - (eA(t,x)e+\epsilon) \partial_{zz}U - \nabla \cdot (A(t,x)e \partial_{z} U)-\partial_{z}(e A(t,x)\nabla U)\\
+q(t,x)\cdot\nabla U+q(t,x)\cdot e\partial_{z}U+c\partial_{z}U+(\beta-\mu) U=0.\\ \end{array}\]
Multiplying this equation by $U$ and integrating by parts over $\Sigma_{a}$ gives:
\[\begin{array}{l}
\frac{1}{2}\big( \int_{(-a,a)\times C}U^{2}(z,T,x)dzdx-\int_{(-a,a)\times C}U^{2}(z,0,x)dzdx\big)\\
=\int_{(-a,a)\times C}\big( -(e\partial_{z}U+\nabla U)A(e\partial_{z}U+\nabla U)-\epsilon\partial_{z}U A\partial_{z}U+\frac{1}{2}\nabla\cdot q U^{2} +(\mu-\beta) U^{2}\big).\\ \end{array}\]
Thus, the choice of $\beta$ yields that:
\[\int_{(-a,a)\times C}U^{2}(z,T,x)dzdx\leq \int_{(-a,a)\times C}U^{2}(z,0,x)dzdx,\]
and thus:
\[\|u_{2}(T)-u_{1}(T)\|_{\mathcal{L}}\leq e^{-\mu T}\|u_{2}(0)-u_{1}(0)\|_{\mathcal{L}}.\]
This means that $G$ is a contraction. The Picard fixed point theorem yields that it admits a unique fixed point $u$, that is, a space-time periodic function $u$ such that $L_{\epsilon}u+\beta u=g$. The Schauder parabolic estimates give the pointwise boundary conditions:
\[u(-a,t,x)=u(a,t,x)=0 \ \hbox{for all} \ (t,x)\in\R\times\R^{N}.\]

Since $\zeta_{\epsilon}\geq\theta_{\epsilon}$ in $\Sigma_{a}$, one can now carry out some iteration procedure. As $g$ is locally  Lipshitz continuous, one can take $\beta$ large enough so that for all $(t,x)\in\R\times\R^{N}$, $s\mapsto f(t,x,s)+\beta s$ is increasing. We define the sequence $(\phi_{n})_{n}$ by:
\begin{equation} \left\{ \begin{array}{rcl} \phi_{0}&=& \zeta_{\epsilon},\\
L_{\epsilon}\phi_{n+1}+\beta\phi_{n+1}&=& f(t,x,\phi_{n})+\beta \phi_{n},\\
\phi_n(-a,t,x)&=& \theta_\epsilon(-a,t,x),\\
\phi_n(a,t,x)&=& \zeta_\epsilon(a,t,x).\\
\end{array}\right. \end{equation}
One can easily prove that the sequence $(\phi_{n})_{n}$ is nonincreasing with respect to $n$ and that for all $n$, one has $\theta_{\epsilon}\leq \phi_{n}\leq\zeta_{\epsilon}$. Thus this sequence converges to some $\phi$ which is a solution of $(\ref{regfinite})$.

It is only left to prove that for all $n$, $\phi_{n}$ is nondecreasing with respect to $z$. We prove this property by iteration. It is clear for $n=0$. Assume that $\phi_{n}$ is nondecreasing with respect to $z$, set $\Sigma_{a}^{\lambda}= (-a,-a+\lambda)\times (0,T)\times C$ and 
\[\phi_{n}^{\lambda}(z,t,x)=\phi_{n}(z+\lambda,t,x),\]
\[\phi_{n+1}^{\lambda}(z,t,x)=\phi_{n+1}(z+\lambda,t,x),\]
for all $\lambda\in (0,2a)$.
In order to show that $\phi_{n+1}$ is nondecreasing in $z$ in $\overline{\Sigma_{a}}$, one only has to show that $\phi_{n+1}\leq\phi_{n+1}^{\lambda}$ in $\overline{\Sigma_{a}^{\lambda}}$ for $\lambda>0$ sufficiently small. 

One can remark that
\[\begin{array}{rcl}(L_{\epsilon}+\beta)(\phi_{n+1}^{\lambda}-\phi_{n+1})&=&f(t,x,\phi_{n}^{\lambda})+\beta\phi_{n}^{\lambda}-f(t,x,\phi_{n})-\beta\phi_{n}\geq 0\\ \end{array}\]
since $s\mapsto f(t,x,s)+\beta s$ is increasing and $\phi_{n}^{\lambda}\geq \phi_{n}$. 

In the other hand, there exists some $a_{0}$ such that for all $a>a_{0}$, the function $z\mapsto \theta_{\epsilon}(z,t,x)$ is increasing in $z\in(-\infty,-a_{0})$ for all $(t,x)\in\R\times\R^{N}$. Fix $a>a_{0}$ and $\lambda\in (0,a-a_{0})$. One has $\phi_{n+1}(-a+\lambda,t,x)-\theta_{\epsilon}(-a,t,x) \geq 0$ for all $(t,x)\in\R\times\R^{N}$. Similarly, as $\zeta_{\epsilon}$ is nondecreasing, one has $\zeta_{\epsilon}(a,t,x)-\phi_{n+1}(a-\lambda,t,x) \geq 0$. This finally gives
\[\phi_{n+1}^{\lambda}(-a,t,x)-\phi_{n+1}(-a,t,x)=\phi_{n+1}(-a+\lambda,t,x)-\theta_{\epsilon}(-a,t,x) \geq 0,\]
\[\phi_{n+1}^{\lambda}(a-\lambda,t,x)-\phi_{n+1}(a-\lambda,t,x)=\zeta_{\epsilon}(a,t,x)-\phi_{n+1}(a-\lambda,t,x) \geq 0.\]
Thus, as $\beta>0$, the strong maximum principle yields that $\phi_{n+1}^{\lambda}\geq \phi_{n+1}$ in $\overline{\Sigma_{a}^{\lambda}}$. Thus $\phi_{n+1}$ is nondecreasing for all $n$. $\Box$


\subsection{Passage to the limit in infinite cylinders}

Let $a_{n} \rightarrow +\infty$ be any sequence that goes to infinity. From standard parabolic estimates and Sobolev's injections, the functions $\phi_{\epsilon,a_{n}}$ converge (up to the extraction of a subsequence) in $\mathcal{C}_{loc}^{2+\beta,1+\beta/2,2+\beta}(\mathbb{R}\times\mathbb{R}\times\mathbb{R}^{N})$, for all $0\leq \beta<\delta$, to a function $\phi_{\epsilon}$ that satisfies:

\begin{equation} \label{reginfinite} \left\{\begin{array}{l}
 L_{\epsilon}\phi_{\epsilon}=f(t,x,\phi_{\epsilon}) \ \hbox{in} \ \mathbb{R}\times\mathbb{R}\times\mathbb{R}^{N}, \\
 \phi_{\epsilon} \ \hbox{is periodic in t and x},\\
 \phi_{\epsilon} \ \hbox{is nondecreasing in z},\\
\end{array} \right. \end{equation}

with $ \theta_{\epsilon}(z,t,x) \leq \phi_{\epsilon}(z,t,x) \leq \zeta_{\epsilon}(z,t,x)$ for all $(z,t,x)\in \mathbb{R}\times\mathbb{R}\times\mathbb{R}^{N}$.

\begin{prop}
The function $\phi_{\epsilon}$ has the following asymptotic behaviours:
\begin{equation}\left\{
\begin{array}{l}
 \phi_{\epsilon}(z,t,x) \rightarrow 0 \ \hbox{as} \ z \rightarrow -\infty\\
 \phi_{\epsilon}(z,t,x) - p(t,x) \rightarrow 0 \ \hbox{as} \ z \rightarrow +\infty\\
\end{array} \right.
\end{equation}
in $\mathcal{C}_{loc}^{1,2}(\mathbb{R}\times\mathbb{R}^{N})$ in the both cases.
\end{prop}

\noindent {\bf Proof.}
From standard parabolic estimates, from the monotonicity of $\phi_{\epsilon}$ in $z$ and from the periodicity in $t$ and $x$, it follows that:
\[ \phi_{\epsilon}(z,t,x) \rightarrow \phi_{\pm}(t,x) \ \hbox{in}\ \mathcal{C}_{loc}^{1,2}(\mathbb{R}\times\mathbb{R}^{N}) \ \hbox{as} \ z \rightarrow \pm \infty, \]
where each function $\phi_{\pm}$ satisfies:
\[ \left\{ \begin{array}{l}
 \partial_{t}\phi_{\pm} -\nabla \cdot(A(t,x)\nabla \phi_{\pm})+q(t,x)\cdot\nabla\phi_{\pm}=f(t,x,\phi_{\pm}),\\
 \phi_{\pm} \ \hbox{is periodic in t and x},\\
 0\leq \phi_{\pm} \leq p.\\
\end{array}
\right. \]

Hypothesis $\ref{uniquenesshyp}$ yields that either $\phi_{\pm} \equiv 0$ or $\phi_{\pm} \equiv p$. Because of the monotonicity of $\phi_{\epsilon}$ in $z$, the inequalities \[0 \leq \phi_{-} \leq \phi_{\epsilon}(z,.,.) \leq \ \phi_{+} \leq p\] hold for all $z$.

If $\phi_{+} \equiv 0$, then $\phi_{\epsilon}(z,t,x) \equiv 0$ for all $(z,t,x) \in \mathbb{R}\times\mathbb{R}\times\mathbb{R}^{N}$. This contradicts the inequality $\phi_{\epsilon} \geq \theta_{\epsilon}$ since $\theta_\epsilon$ is not uniformly nonpositive.
This shows that $\phi_{\epsilon}(z,t,x) \rightarrow p(t,x)$ as $z \rightarrow +\infty$.

Similarly, if $\phi_{-} \equiv p$ then $\phi_{\epsilon} \equiv p$. This contradicts the inequality $\phi_{\epsilon}(z,t,x)\leq \zeta_\epsilon $ when $z$ goes to $-\infty$ since $p$ is periodic and positive. This shows that $\phi_{\epsilon}(z,t,x) \rightarrow 0$ as $z \rightarrow -\infty$. $\Box$


\subsection{Removal of the regularization}

In this section, our aim is to let $\epsilon\rightarrow 0$. The following result is true even if $f$ is not of KPP type. 

\begin{prop} \label{removeregpos}
Assume that $p$ is a space-time periodic positive solution of equation $(\ref{eqprinc})$.
Consider $(\phi_{\epsilon})_{\epsilon \in\mathcal{E}}$ a family of solutions of equation $(\ref{reginfinite})$, where $\mathcal{E}$ is a subset of $\R^{+}$, such that for all $\epsilon \in\mathcal{E}$, $\phi_{\epsilon}$ is nondecreasing with respect to $z$.
Then, the family $(\partial_{z}\phi_{\epsilon})_{\epsilon \in\mathcal{E}}$ is uniformly bounded in $L^{1}_{loc}(\R\times\R\times\R^{N})$ and the families $(\nabla \phi_{\epsilon}+e\partial_{z}\phi_{\epsilon})_{\epsilon \in\mathcal{E}}$ and $(\partial_{t}\phi_{\epsilon}+c\partial_{z}\phi_{\epsilon})_{\epsilon \in\mathcal{E}}$ are uniformly bounded in $L^{2}_{loc}(\R\times\R\times\R^{N})$. Furthermore, these bounds are locally uniform with respect to $c$.
\end{prop}

\noindent {\bf Proof.}
First of all, as $\phi_{\epsilon}$ is nondecreasing in $z$, for all $R>0, \epsilon \in\mathcal{E}$, we have:
\begin{equation} \begin{array}{l} 
\int_{(-R,R)\times (0,T) \times C} |\partial_{z}\phi_{\epsilon}|dzdtdx =\int_{(-R,R)\times (0,T) \times C}\partial_{z}\phi_{\epsilon}dzdtdx \\
\\
= \int_{(0,T)\times C}\phi_{\epsilon}(R,t,x)dtdx -\int_{(0,T)\times C}\phi_{\epsilon}(-R,t,x)dtdx \rightarrow \int_{(0,T)\times C} p(t,x) dtdx \ \hbox{as} \ R\rightarrow +\infty \\
\end{array}
\end{equation}
which proves that $\partial_{z}\phi_{\epsilon}$ is uniformly bounded in $L^{1}_{loc}(\R\times\R\times\R^N)$.  Multiplying equation $(\ref{reginfinite})$ by $\phi_{\epsilon}$ and integrating, we get for all $R>0, \epsilon\in\mathcal{E}$:
\begin{equation} \begin{array}{l}
\int_{(-R,R)\times (0,T)\times C}(\nabla \phi_{\epsilon}+e\partial_{z}\phi_{\epsilon})A(t,x)(\nabla \phi_{\epsilon}+e\partial_{z}\phi_{\epsilon})dzdtdx+\epsilon\int_{(-R,R)\times (0,T)\times C}(\partial_{z}\phi_{\epsilon})^{2}dzdtdx\\
\\
+\int_{(-R,R)\times (0,T) \times C}(c+q\cdot e)\partial_{z}(\frac{\phi_{\epsilon}^{2}}{2})dzdtdx-\frac{1}{2}\int_{(-R,R)\times (0,T)\times C}(\nabla\cdot q)\phi_{\epsilon}^{2}dzdtdx\\
\\
= \int_{(-R,R)\times (0,T)\times C}\phi_{\epsilon}f(t,x,\phi_{\epsilon})dzdtdx \\
\end{array} \end{equation}

Using the ellipticity property of the matrix $A$ and the inequality $0\leq \phi_{\epsilon}\leq p$, we get:
\begin{equation} \begin{array}{rcl}
&&\gamma\int_{(-R,R)\times (0,T)\times C}|\nabla \phi_{\epsilon}+e\partial_{z}\phi_{\epsilon}|^{2}dzdtdx\\
&&\\
&\leq& \int_{(-R,R)\times (0,T)\times C}(\phi_{\epsilon}f(t,x,\phi_{\epsilon})+\frac{\nabla\cdot q}{2}\phi_{\epsilon}^{2})dzdtdx\\
&&+\int_{(0,T) \times C}(c+q\cdot e)(\frac{\phi_{\epsilon}^{2}}{2}(-R,t,x) 
-\frac{\phi_{\epsilon}^{2}}{2}(R,t,x))dtdx \\
\\
&\leq& (M+\frac{1}{2}\|\nabla\cdot q\|_{\infty})\int_{(-R,R)\times (0,T)\times C}|\phi_{\epsilon}|^{2}dzdtdx+\int_{(0,T) \times C}|c+q\cdot e|p^{2}(t,x)dtdx \\
\end{array} \end{equation}
where $\eta(t,x)=\sup_{0<s<p(t,x)}\displaystyle\frac{f(t,x,s)}{s}$ and $M=\sup_{(t,x)\in\R\times\R^N}|\eta(t,x)|<\infty$ since $f$ is of class $\mathcal{C}^1$. Finally, this gives:
\[\|\nabla\phi_{\epsilon}+e\partial_{z}\phi_{\epsilon}\|_{L^{2}((-R,R)\times (0,T)\times C)}\leq \frac{1}{\sqrt{\gamma}} \big(2R(\|\eta\|_{\infty}+\frac{1}{2}\|\nabla\cdot q\|_{\infty})+|c|+\|q\|_{\infty}\big)^{1/2}\|p\|_{L^{2}((0,T)\times C)}.\]
It follows that $\nabla \phi_{\epsilon}+e\partial_{z}\phi_{\epsilon}$ is uniformly bounded in $L^{2}((-R,R)\times (0,T)\times C)$ for all positive $R$.

Similarly, multiplying equation $(\ref{reginfinite})$ by $\partial_{t}\phi_{\epsilon}+c\partial_{z}\phi_{\epsilon}$ and integrating, we get for all $R>0, \epsilon\in\mathcal{E}$: 
\begin{equation} \begin{array}{rcl}
&&\int_{(-R,R)\times(0,T)\times C}|\partial_{t}\phi_{\epsilon}+c\partial_{z}\phi_{\epsilon}|^{2}dzdtdx\\
&&\\
&=&-\int_{(-R,R)\times(0,T)\times C}q\cdot(\nabla\phi_{\epsilon}+e\partial_{z}\phi_{\epsilon})(\partial_{t}\phi_{\epsilon}+c\partial_{z}\phi_{\epsilon})dzdtdx\\
&&\\
&&+\int_{(-R,R)\times(0,T)\times C}\{(e\partial_{z}+\nabla)(A(t,x)e\partial_{z}\phi_{\epsilon}+A(t,x)\nabla\phi_{\epsilon})\}(\partial_{t}\phi_{\epsilon}+c\partial_{z}\phi_{\epsilon})dzdtdx\\
&&\\
&&+\epsilon\int_{(-R,R)\times(0,T)\times C}\partial_{zz}\phi_{\epsilon}(\partial_{t}\phi_{\epsilon}+c\partial_{z}\phi_{\epsilon})dzdtdx\\
&&\\
&&+\int_{(-R,R)\times(0,T)\times C}f(t,x,\phi_{\epsilon})(\partial_{t}\phi_{\epsilon}+c\partial_{z}\phi_{\epsilon})dzdtdx\\
&&\\
\end{array}
\end{equation}

\begin{equation}
\begin{array}{rcl}
&\leq&\|q\|_{L^{\infty}}\|\nabla\phi_{\epsilon}+e\partial_{z}\phi_{\epsilon}\|_{L^{2}((-R,R)\times (0,T)\times C)} \|\partial_{t}\phi_{\epsilon}+c\partial_{z}\phi_{\epsilon}\|_{L^{2}((-R,R)\times (0,T)\times C)}\\
&&\\
&&+\int_{(0,T)\times C}(e\partial_{z}\phi_{\epsilon}+\nabla\phi_{\epsilon})A e(\partial_{t}\phi_{\epsilon}+c\partial_{z}\phi_{\epsilon})(R,t,x)dtdx\\
&&\\
&&-\int_{(0,T)\times C}(e\partial_{z}\phi_{\epsilon}+\nabla\phi_{\epsilon})A e(\partial_{t}\phi_{\epsilon}+c\partial_{z}\phi_{\epsilon})(-R,t,x)dtdx \\
&&\\
&&-\int_{(-R,R)\times(0,T)\times C}(e\partial_{z}\phi_{\epsilon}+\nabla\phi_{\epsilon})A(t,x)(\partial_{t}+c\partial_{z})(e\partial_{z}\phi_{\epsilon}+\nabla\phi_{\epsilon})dzdtdx\\
&&\\
&&+\epsilon\int_{(0,T)\times C}\big(\partial_{z}\phi_{\epsilon}\partial_{t}\phi_{\epsilon}(R,t,x)-\partial_{z}\phi_{\epsilon}\partial_{t}\phi_{\epsilon}(-R,t,x)\big)dtdx\\
&&\\
&&+\frac{c\epsilon}{2}\int_{(0,T)\times C}\big( (\partial_{z}\phi_{\epsilon})^{2}(R,t,x)-(\partial_{z}\phi_{\epsilon})^{2}(-R,t,x)\big)dtdx\\
&&\\
&&+\big(\int_{(-R,R)\times(0,T)\times
C}f(t,x,\phi_{\epsilon})^{2}dzdtdx)^{1/2}\|\partial_{t}\phi_{\epsilon}+c\partial_{z}\phi_{\epsilon}\|_{L^{2}((-R,R)\times (0,T)\times C)}.\\
&&\\
\end{array}
\end{equation}

\begin{equation}
\begin{array}{rcl}
&\leq&\|q\|_{L^{\infty}}\|\nabla\phi_{\epsilon}+e\partial_{z}\phi_{\epsilon}\|_{L^{2}((-R,R)\times (0,T)\times C)}\|\partial_{t}\phi_{\epsilon}+c\partial_{z}\phi_{\epsilon}\|_{L^{2}((-R,R)\times (0,T)\times C)}\\
&&\\
&&+\frac{1}{2}\int_{(-R,R)\times(0,T)\times C}(e\partial_{z}\phi_{\epsilon}+\nabla\phi_{\epsilon})\partial_{t}A(e\partial_{z}\phi_{\epsilon}+\nabla\phi_{\epsilon})dzdtdx\\
&&\\
&&\sqrt{2R}\|\eta\|_{\infty}\|p\|_{L^{2}((0,T)\times C)}\|\partial_{t}\phi_{\epsilon}+c\partial_{z}\phi_{\epsilon}\|_{L^{2}((-R,R)\times (0,T)\times C)},\\
&&\\
\end{array}
\end{equation}

\begin{equation}
\begin{array}{rcl}
&\leq&\big(\|q\|_{L^{\infty}}\|\nabla\phi_{\epsilon}+e\partial_{z}\phi_{\epsilon}\|_{L^{2}}+\sqrt{2R}\|\eta\|_{\infty}\|p\|_{L^{2}((0,T)\times C)}\big)\|\partial_{t}\phi_{\epsilon}+c\partial_{z}\phi_{\epsilon}\|_{L^{2}((-R,R)\times (0,T)\times C)}\\
&&\\
&&+\frac{1}{2}\|\partial_{t}A\|_{\infty}\|e\partial_{z}\phi_{\epsilon}+\nabla\phi_{\epsilon}\|_{L^{2}((-R,R)\times\R\times\R^{N})}^{2}\\
\end{array}\end{equation}

Thus, as the previous estimates yield that $\nabla\phi_{\epsilon}+e\partial_{z}\phi_{\epsilon}$ is uniformly bounded in $L^{2}((-R,R)\times (0,T)\times C)$ for all positive $R$, this computation proves that $\partial_{t}\phi_{\epsilon}+c\partial_{z}\phi_{\epsilon}$ is uniformly bounded in $L^{2}((-R,R)\times\R\times\R^{N})$ for all positive $R$. $\Box$

\bigskip

We now apply this theorem and pass to the limit $\epsilon\rightarrow 0$.

\noindent {\bf Proof of Theorem \ref{thmexistence} in the case $c>c^*_e$.}
As $L^{2}_{loc}(\R\times\R\times\R^{N})$ is embedded in $L^{1}_{loc}(\R\times\R\times\R^{N})$, Proposition $\ref{removeregpos}$ applied to $\mathcal{E}=\R^{+*}$ and the diagonal extraction process yield that there exists a sequence $(\epsilon_{n})$ that converges to $0$ and a limit function $\phi$ such that:
\[ \phi_{\epsilon_{n}} \rightarrow \phi \ \hbox{in} \ L^{1}_{loc}(\R\times\R\times\R^{N}) \ \hbox{and almost everywhere}.\]
This convergence yields that $\phi$ solves the degenerate equation $(\ref{eqtf})$ in the sense of distributions and is nondecreasing almost everywhere in $z$, that is, for almost every $(z_{1},z_{2})\in\R^{2}$ and $(t,x)\in\R\times\R^{N}$, one has $\phi(z_{1},t,x)\leq \phi (z_{2},t,x)$. As $\phi_{\epsilon}\leq p$ for all $\epsilon>0$, one has $\phi\leq p$ almost everywhere and thus $\phi\in L^{\infty}(\R\times\R\times\R^{N})$. 
Furthermore, one has
\begin{equation} \label{normalization}\theta_{0}(z,t,x)\leq \phi (z,t,x)\leq \psi_{\mu,\lambda_{c}^{0}(\mu)}(t,x)e^{\lambda_{c}^{0}(\mu)z}\end{equation}
for almost every $(z,t,x)\in \R\times\R\times\R^{N}$. This gives $(\ref{expdec})$.
It is only left to prove that the asymptotic conditions are satisfied. 

\begin{prop} \label{asymptotics2}
Consider a solution $\phi\in W^{1,1}_{loc}(\R\times\R\times\R^{N})$ of equation $(\ref{eqtf})$ in the sense of distributions, which is nondecreasing almost everywhere in $z$. Assume that this function is not uniformly equal to $0$ or to $p$. Then it satisfies the following limits as $z \rightarrow \pm \infty$:
\[ \phi(z,t,x) \rightarrow 0 \ \hbox{as} \ z\rightarrow -\infty \ \hbox{and} \ \phi(z,t,x)-p(t,x) \rightarrow 0 \ \hbox{as} \ z\rightarrow +\infty \]
in $L^{\infty}_{loc}(\R\times\R^{N})$.
\end{prop}

This proposition concludes the proof of Theorem $\ref{thmexistence}$ in the case $c>c^{*}_e(A,q,\mu)$ since our function $\phi$ is not uniformly equal to $0$ or $p$ thanks to $(\ref{normalization})$. We will conclude this proof in the case $c=c^{*}_e(A,q,\mu)$ later using Proposition $\ref{vitferme}$.

\bigskip

\noindent {\bf Proof.}
As $\phi$ is nondecreasing almost everywhere in $z$ and bounded, one can assume, up to some change of this function on a null-measure set, that $\phi$ is nondecreasing. Then there exist two periodic functions $\phi_{\pm}$ such that: 
\[ \phi(z,t,x) \rightarrow \phi_{\pm}(t,x) \ \hbox{as} \ z\rightarrow \pm \infty,\]
for almost every $(t,x)\in\R\times\R^{N}$. 
It is left to prove that $\phi_{-}\equiv 0$ and $\phi_{+}\equiv p$.

Take $h$ a smooth function periodic in $t$ and $x$. Take $\xi_{0}\in \mathcal{C}^{\infty}(\R)$ a nonnegative bounded function that satisfies:
\[\xi_{0}(z)=0 \ \hbox{if} \ |z|\leq 1 \hbox{ and } \int_{\R}\xi_{0}(z)dz = 1,\]
and for all $n\in\N$, set $\xi_{n}(z)=\xi_{0}(z-n)$. 

As $\phi$ is a weak solution of $(\ref{eqtf})$, multiplying $(\ref{eqtf})$ by $\xi_{n}(z)h(t,x)$ and integrating over $\R\times (0,T)\times C$ gives:
\begin{equation} \label{cvdeg}\begin{array}{l}
\int_{\R\times (0,T)\times C}f(t,x,\phi)h(t,x)\xi_{n}(z)dzdtdx\\
\\
=\int_{\R\times (0,T)\times C} (\partial_{t}\phi-\nabla \cdot(A\nabla\phi)-2eA\nabla\partial_{z}\phi -eAe\partial_{zz}\phi)h(t,x)\xi_n(z)dzdtdx\\
+\int_{\R\times (0,T)\times C}q\cdot\nabla\phi+c\partial_{z}\phi+q\cdot e\partial_{z}\phi)h(t,x)\xi_n(z)dzdtdx\\
\\
=\int_{(0,T)\times C}(-\partial_{t}h-\nabla\cdot (A\nabla h)-\nabla\cdot (qh))(\int_{\R}\phi(z,t,x)\xi_{n}(z)dz)dtdx\\
\\
\ +\int_{\R\times (0,T)\times C} (2eA\nabla\phi h(t,x) \xi_{n}'(z)-eAe\phi\xi_{n}''(z)h(t,x)-q\cdot e \phi h(t,x)\xi_{n}'(z))dzdtdx.\\
\end{array}\end{equation}

One can compute:
\begin{equation} \begin{array}{rcl}
\int_{\R\times (0,T)\times C} eAe\phi h \xi_{n}''dzdtdx &=& \int_{(0,T)\times C}eA(t,x)e h(t,x) (\int_{\R}\phi(z,t,x)\xi_{n}''(z)dz)dtdx\\
&=& \int_{(0,T)\times C}eA(t,x)e h(t,x) (\int_{-1}^{1}\phi(z+n,t,x)\xi_{0}''(z)dz)dtdx\\
&\rightarrow & \int_{(0,T)\times C}eA(t,x)e h(t,x) (\int_{-1}^{1}\phi_{+}(t,x)\xi_{0}''(z)dz)dtdx
\end{array} \end{equation}
as $n\rightarrow +\infty$, and thus:
\[\int_{\R\times (0,T)\times C} eAe\phi h \xi_{n}''dzdtdx\rightarrow \int_{(0,T)\times C}eA(t,x)e h(t,x)\phi_{+}(t,x) (\int_{\R}\xi_{0}''(z)dz)dtdx =0.\]

Computing each term of equation $(\ref{cvdeg})$ in a similar way, one gets:
\begin{equation}\begin{array}{l}
\int_{\R\times (0,T)\times C} (\partial_{t}\phi-\nabla \cdot(A\nabla\phi)-2eA\nabla\partial_{z}\phi -eAe\partial_{zz}\phi+q\cdot\nabla\phi+c\partial_{z}\phi+q\cdot e\partial_{z}\phi)dzdtdx\\
\rightarrow \int_{(0,T)\times C}(-\partial_t h-\nabla\cdot (A\nabla h)-\nabla\cdot (qh))\phi_{+}dtdx \ \hbox{as} \ n\rightarrow +\infty\\
\end{array}\end{equation}
and
\begin{equation} \begin{array}{l}
\int_{\R\times (0,T)\times C}f(t,x,\phi)h(t,x)\xi_{n}(z)dzdtdx
\rightarrow \int_{(0,T)\times C}f(t,x,\phi_{+})h(t,x)dtdx\\
\end{array}\end{equation}
as $n\rightarrow +\infty$. This yields that $\phi_{+}$ is a weak solution of the equation:
\[\partial_{t}\phi-\nabla\cdot (A\nabla\phi)+q\cdot \nabla\phi=f(t,x,\phi).\]

The regularity theorem for parabolic equations yields that this is a strong periodic solution of this equation. 
Hypothesis $\ref{uniquenesshyp}$ yields that there only exist two periodic nonnegative solutions of this equation: $0$ and $p$. 
As $\phi$ is nondecreasing almost everywhere, one has: $0\leq\phi_{-}(t,x)\leq\phi(z,t,x)\leq\phi_{+}(t,x)\leq p(t,x)$ for almost every $(z,t,x)\in\R\times\R\times\R^{N}$.  

Assume that $\phi_{+}\equiv 0$, then one has $\phi(z,.,.)\equiv 0$ for almost every $z$. This contradicts the estimate $(\ref{normalization})$. Thus, $\phi_{+}\equiv p$. Similarly, one can prove that $\phi_{-}\equiv 0$.

The Dini's lemma and the periodicity give that the previous convergence as $z\rightarrow+\infty$ is uniform in $(t,x)\in\R\times\R^{N}$.
Similarly, one can prove that $\phi(z,t,x)  \rightarrow 0$ as $z\rightarrow -\infty$ uniformly in $(t,x)\in\R\times\R^N$. 
$\Box$


\section{The general case} 

\subsection{Proof of the existence result}

\noindent {\bf Proof of Theorem \ref{thmexistencepos} in the case $c>c^*_e$.}
We fix $c>c^{*}_{e}(A,q,\eta)$ and $\epsilon_{0}$ such that for all $0<\epsilon<\epsilon_{0}$, one has $c>c^*_{\epsilon}(\eta)$.
Hence $\lambda_{c}^{\epsilon}(\eta)$ is well-defined by Proposition \ref{deflambdac} and one can set \[\zeta_{\epsilon}(z,t,x)=\inf \{ p(t,x), \psi_{\eta,\lambda_{c}^{\epsilon}(\eta)}(t,x)e^{\lambda_{c}^{\epsilon}(\eta)z}\},\] 
where $\psi_{\eta,\lambda}$ is the unique space-time periodic principal eigenfunction defined by (\ref{eigen}) but with the zero order term $\eta(t,x)=\sup_{0<s<p(t,x)}f(t,x,s)/s$ and normalized by $\|\psi_{\eta,\lambda}\|_\infty=1$. 

As $\eta\leq\mu$ and $\eta\not\equiv \mu$, one has $k_{\lambda e}(A,q,\mu)< k_{\lambda e}(A,q,\eta)$ and thus $c^{*}_{\epsilon}(\mu)\leq c^{*}_{\epsilon}(\eta)<c$. Hence the function $\theta_{\epsilon}$ that was used in the previous section (see Lemma $\ref{thetalem}$ for the definition) is still well-defined.
Up to some translation, we assume that \[\max_{z\in\R}\min_{(t,x)\in\R\times\R^{N}}\theta_{\epsilon}(z,t,x)=\min_{(t,x)\in\R\times\R^{N}}\theta_{\epsilon}(0,t,x).\] Thus $z\mapsto \theta_{\epsilon}(z,t,x)$ is increasing for all $(t,x)\in\R\times\R^{N}$ over $z\in\R^{-}$.

The function $\zeta_{\epsilon}$ decreases to $0$ with the rate $\lambda_{c}^{\epsilon}(\eta)$ and the function $\theta_{\epsilon}$ decreases with the rate $\lambda_{c}^{\epsilon}(\mu)$. As $\mu\leq\eta$, it is possible to prove that $\lambda_{c}^{\epsilon}(\mu)\leq\lambda_{c}^{\epsilon}(\eta)$ and thus one cannot expect to get $\theta_{\epsilon}\leq\zeta_{\epsilon}$ in $\R$. Anyway, it is still possible to get such a comparison on finite intervals $(-a,a)$. 

We thus investigate the approximated problem:
\begin{equation} \label{regfinitepos}
\left\{ \begin{array}{rcl}
L_{\epsilon} \phi &=& f(t,x,\phi), \\
\phi \ \hbox{is periodic} \ &\hbox{in}& \ t,x, \\
\phi(-a,t,x) &=& \theta_{\epsilon}(-a+m_{a}(\tau),t,x) \ \hbox{for all} \ (t,x) \in \mathbb{R} \times \mathbb{R}^{N},\\
\phi(a,t,x) &=& \zeta_{\epsilon}(a+\tau,t,x) \ \hbox{for all} \ (t,x) \in \mathbb{R} \times \mathbb{R}^{N}, \\
\end{array}
\right.
\end{equation}
where $m_{a}(\tau)$ is defined by
\[m_a(\tau)=\min \big\{0, \frac{\lambda_{\epsilon}^{c}(\eta)}{\lambda_{\epsilon}^{c}(\mu)}(\tau-a)+a\big\}.\]
As $\eta\geq\mu$, one has $\lambda_{\epsilon}^{c}(\eta)\geq\lambda_{\epsilon}^{c}(\mu)$. Thus $\zeta_{\epsilon}$ increases faster than $\theta_{\epsilon}$. As $m_a(\tau)$ has been chosen so that for all $(t,x)\in\R\times\R^{N}$:
\[\zeta_{\epsilon}(-a+\tau,t,x)\geq \theta_{\epsilon}(-a+m_a(\tau),t,x),\]
one finally gets 
\[\zeta_{\epsilon}(z+\tau,t,x)\geq \theta_{\epsilon}(z+m_a(\tau),t,x) \ \hbox{for all} \ (z,t,x)\in (-a,+\infty) \times\R\times\R^{N}.\]

Thus one can use the same method as in the proof of Lemma $\ref{existencerefinite}$ with the subsolution $(z,t,x)\mapsto\theta_{\epsilon}(z+m_{a}(\tau),t,x)$ and the supersolution $(z,t,x)\mapsto\zeta_{\epsilon}(z+\tau,t,x)$ to prove the existence of a solution $\phi_{\epsilon,a}^{\tau}$ of equation $(\ref{regfinitepos})$
that satisfies \[\theta_{\epsilon}(z+m_{a}(\tau),t,x)\leq \phi_{\epsilon,a}^{\tau}(z,t,x)\leq \zeta_{\epsilon}(z+\tau,t,x)\] for all $(z,t,x)\in\Sigma_{a}$. Moreover, as $(z,t,x)\mapsto\theta_{\epsilon}(z+m_{a}(\tau),t,x)$ is increasing in the neighborhood of $-a<0$ since $\theta_{\epsilon}$ is increasing with respect to $z\in\R^{-}$ and $m_{a}(\tau)\leq 0$, it is possible to choose some $\phi_{\epsilon,a}^{\tau}$ which is nondecreasing with respect to $z$.

As the boundary conditions are continuous with respect to $\tau$, the Schauder interior estimates give that $\phi_{\epsilon,a}^{\tau}$ is continuous with respect to $\tau$. Similarly, the boundary conditions are nondecreasing with respect to $\tau$ since $\theta_{\epsilon}$ is nondecreasing with respect to $z\in\R^{-}$ and $m_{a}(\tau)\leq 0$. Thus, using a sliding method as in the proof of the monotonicity in Lemma $\ref{existencerefinite}$, one gets that $\tau\mapsto\phi_{\epsilon,a}^{\tau}(z,t,x)$ is nondecreasing for all $(z,t,x)\in (-a,a)\times\R\times\R^{N}$. As $\zeta_{\epsilon}(-\infty,t,x)=0$, the function $\phi_ {\epsilon,a}^{\tau}$ uniformly converges to $0$ as $\tau\rightarrow -\infty$ in $(-a,a)\times\R\times\R^{N}$. For all $\tau\geq 0$ and $z\geq 0$, one has
\[\phi_{\epsilon,a}^{\tau}(z,t,x)\geq \phi_{\epsilon,a}^{\tau}(0,t,x)\geq \phi_{\epsilon,a}^{0}(0,t,x)\geq \theta_{\epsilon}(0,t,x) \ \hbox{for all} \ (t,x)\in\R\times\R^{N}.\]
Set $\theta^{-}=\min_{(t,x)\in\R\times\R^{N}, \epsilon>0}\theta_{\epsilon}(0,t,x)$. One can fix some $\tau=\tau_{\epsilon,a}$ such that:
\[\frac{1}{T|C|}\int_{(0,1)\times (0,T)\times C} \phi_{\epsilon,a}^{\tau}(z,t,x)dzdtdx=\frac{\theta^{-}}{2}.\]

As $a\rightarrow +\infty$, one may assume, up to extraction, that $\phi_{\epsilon,a}^{\tau_{\epsilon,a}}$ converges to some function $\phi_{\epsilon}$ in $\mathcal{C}^{2,1,2}_{loc}(\R\times\R\times\R^{N})$. This function is periodic in $(t,x)$, nondecreasing in $z$ and satisfies:
\[\left\{ \begin{array}{l} L_{\epsilon}\phi_{\epsilon}=f(t,x,\phi_{\epsilon}),\\
\phi_{\epsilon} \ \hbox{is periodic in} \ (t,x), \\
\phi_{\epsilon} \ \hbox{is increasing in} \ z,\\
\frac{1}{T|C|}\int_{(0,1)\times (0,T)\times C} \phi_{\epsilon}(z,t,x)dzdtdx=\frac{\theta^{-}}{2}.\\ \end{array} \right.\]
Define $\phi^{\pm}(t,x)=\lim_{z\rightarrow \pm \infty}\phi_{\epsilon}(z,t,x)$. Using Hypothesis $\ref{uniquenesshyp}$, one can prove that $\phi^{+}\equiv p$ and $\phi^{-}\equiv 0$. 

All the hypotheses of Proposition $\ref{removeregpos}$ are now satisfied and thus one can assume, up to extraction, that $\phi_{\epsilon}$ converges to some function $\phi$ in $L^{1}_{loc}(\R\times\R\times\R^{N})$ such that:
\begin{equation} 
\left\{ \begin{array}{l}
\partial_{t} \phi - \nabla \cdot (A(t,x) \nabla \phi) - eA(t,x)e \partial_{zz}\phi - \nabla \cdot (A(t,x)e \partial_{z} \phi)\\
-\partial_{z}(e A(t,x)\nabla \phi) 
+q(t,x)\cdot\nabla\phi+q(t,x)\cdot e\partial_{z}\phi+c\partial_{z}\phi = f(t,x,\phi), \\
\phi \ \hbox{is periodic in} \ (t,x), \\
\phi \ \hbox{is increasing in} \ z,\\
\frac{1}{T|C|}\int_{(0,1)\times (0,T)\times C} \phi(z,t,x)dzdtdx=\frac{\theta^{-}}{2}.\\
\end{array}
\right.
\end{equation}
Proposition $\ref{asymptotics2}$ gives that:
\[\left\{ \begin{array}{rcl}
\phi(z,t,x)&\rightarrow& 0 \ \hbox{as} \ z\rightarrow -\infty,\\
\phi(z,t,x)-p(t,x)&\rightarrow& 0 \ \hbox{as} \ z\rightarrow +\infty.\\ \end{array} \right. \]
Thus $\phi$ is the profile of a pulsating traveling front of speed $c$ and the proof is done for all $c>c^{*}_e(A,q,\eta)$. The proof will be completed later in the case $c=c^{*}_e(A,q,\eta)$ (see Proposition $\ref{vitferme}$). $\Box$


\subsection{The case $\lambda_{1}'=0$}

In this section we prove Theorem $\ref{thmexistencepos}$ and Proposition $\ref{cexlambda1'=0}$. 

\noindent {\bf Proof of Theorem \ref{thmexistencepos}.}
Fix $c>c^{*}(\eta)$. Take $\chi$ some smooth function such that $\chi(s)=1$ if $s\leq 0$, $\chi(s)=0$ if $s\geq 1$ and $\chi(s)>0$ if $s\in (0,1)$. Set $f_{\epsilon}(t,x,s)=f(t,x,s)+\epsilon\chi (s)s$ and $\eta_{\epsilon}(t,x)=\sup_{0<s<1}\frac{f_{\epsilon}(t,x,s)}{s}$. As $(f_{\epsilon})_u'(t,x,0)=f_u' (t,x,0)+\epsilon$, one gets 
\[\lambda_{1}'(A,q,(f_{\epsilon})_u'(t,x,0))=-\epsilon<0.\]

For $\epsilon>0$ small enough, $\eta_{\epsilon}$ is close to $\eta$ and thus $c^{*}(\eta_{\epsilon})$ is close to $c^{*}(\eta)$. Notice that this quantity is finite since, as $k_{0}(A,q,\eta)=\lambda_{1}'(A,q,\eta)=0$, one has:
\[\frac{-k_{\lambda}}{\lambda}=-\frac{k_{\lambda}-k_{0}}{\lambda}\geq -\partial_{\lambda}k_{0}\]
using the concavity of $\lambda\mapsto k_{\lambda}$. This gives $c^{*}(\eta)\geq -\partial_{\lambda}k_{0}(A,q,\eta)$.

Lastly, for all $\epsilon>0$, the function $f_{\epsilon}$ satisfies the hypotheses of Proposition $\ref{reactionhyppos}$. Thus it satisfies Hypothesis $\ref{uniquenesshyp}$.

For $\epsilon$ small enough, one has $c>c^{*}(A,q,\eta_{\epsilon})$ and thus there exists a pulsating traveling front associated with $(A,q,f_{\epsilon})$ using Theorem $\ref{thmexistence}$. Take $\phi_{\epsilon}$ a profile normalized by 
\[\int_{z\in (0,1)}\int_{(0,T)\times C}\phi_{\epsilon}(z,t,x)dzdtdx = \frac{|C|T}{2}.\]
The estimates that were used in the proof of Proposition $\ref{removeregpos}$ were locally uniform in $f$ and then the sequence $(\phi_{\epsilon})$ is uniformly bounded in $W^{1,1}_{loc}(\R\times\R\times\R^{N})$. Thus one can assume, up to extraction, that this sequence converges almost everywhere and in $L^{1}_{loc}(\R\times\R\times\R^{N})$ to a function $\phi$. This function is a weak solution of equation $(\ref{eqtf})$. Using a similar method as in the proof of  Proposition $\ref{asymptotics2}$, we can prove that the good asymptotic behaviours hold when $z \rightarrow \pm\infty$. $\Box$

\bigskip

In order to prove Proposition $\ref{cexlambda1'=0}$, we begin with the following lemma, which is an extension of Proposition 2.13 of \cite{eigenvalue}:
\begin{lem}
For all $(A,q,\mu)$, one has:
\begin{equation} \label{eigenlemmabounded} k_{0}(A,q,\mu)=\inf \{ \lambda \ | \ \exists \phi \in \mathcal{C}^{1,2}(\mathbb{R}\times\R^{N})\cap W^{1,\infty}(\R\times\R^{N}), \phi >0 \ \hbox{and} \ \mathcal{L}\phi \leq \lambda \phi \ \hbox{in} \ \mathbb{R}\times \R^{N} \}.\end{equation}
\end{lem}

\noindent {\bf Proof.}
We forget the dependence in $(A,q,\mu)$ to simplify the notations and we set $\lambda_{1}''$ the rightmember of $(\ref{eigenlemmabounded})$.
Taking $\varphi_{0}$ a periodic principal eigenfunction associated
with $k_{0}$ and using it as a test-function in $(\ref{eigenlemmabounded})$, one immediately gets
$\lambda_{1}''\leq k_{0}$. Next, take $\lambda<k_{0}$ and assume that
there exists a function $\phi\in \mathcal{C}^{1,2}(\R\times
\R^{N})\cap W^{1,\infty}(\R\times \R^{N})$ such that $\phi$ is positive and satisfies $\mathcal{L}\phi \leq
\lambda \phi$. We now search for a contradiction in order to prove
that such a $\lambda$ does not exist and that $\lambda_{1}''\geq k_{0}$.

Set $\gamma=\sup_{\R\times\R^{N}}\frac{\phi}{\varphi_{0}}$, where $\varphi_{0}$ is some space-time periodic eigenfunction associated with $k_{0}$. Then
$0<\gamma<\infty$ and one can define $z=\gamma\varphi_{0}-\phi$.
This function is nonnegative and $\inf z=0$. Set
$\epsilon=(k_{0}-\lambda)\min \varphi_{0} >0$. One has
$(\mathcal{L}-\lambda)(z)\geq\gamma\epsilon$.

Consider a nonnegative function $\theta \in \mathcal{C}^{2}(\R\times\R^{N})$
that satisfies:
\[ \theta (0,0)=0, \lim_{|t|+|x|\rightarrow +\infty} \theta(t,x)=1,
\|\theta\|_{\mathcal{C}^{1,2}}<\infty. \] There exists $\kappa>0$
sufficiently large such that:
\[\forall (s,y)\in\R\times\R^{N}, (\mathcal{L}-\lambda)(\tau_{s,y}\theta) >
-\kappa\gamma\epsilon/2,\] where we denote $\tau_{s,y}\theta =
\theta(.-s,.-y)$.

Since $\inf z = 0$, one can find some
$(t_{0},x_{0})\in\R\times\R^{N}$ such that:
\[z(t_{0},x_{0}) < \min \{
\frac{1}{\kappa},\frac{\gamma\epsilon}{2\|\mu\|_{\infty}} \}\] where
$\|\mu\|_{\infty}=+\infty$ if $\mu\equiv 0$. Since
$\lim_{|t|+|x|\rightarrow +\infty} \theta(t,x)=1$, there exists a positive
constant $R$ such that $\tau_{t_{0},x_{0}}\theta(t,x)/\kappa >
z(t_{0},x_{0})$ if $|t-t_{0}|+|x-x_{0}|\geq R$. Consequently, setting
$\tilde{z}=z+\tau_{t_{0},x_{0}}\theta(t,x)/\kappa$, one finds for all
$|t-t_{0}|+|x-x_{0}|\geq R$, that:
\[\tilde{z}(t,x)\geq \tau_{t_{0},x_{0}}\theta(t,x)/\kappa >
z(t_{0},x_{0})=\tilde{z}(t_{0},x_{0}).\]

Hence, if $\alpha=\inf_{\R\times \R^{N}} \tilde{z}$, this
infimum is reached in $B_{R}(t_{0},x_{0})$. Moreover:
\[\alpha \leq \tilde{z}(t_{0},x_{0}) = z(t_{0},x_{0})
<\frac{\gamma\epsilon}{2\|\mu\|_{\infty}}.\] One can compute:
\[\begin{array}{rcl}(\mathcal{L}-\lambda)(\tilde{z}-\alpha)&=&(\mathcal{L}-\lambda)(z)+\frac{1}{\kappa}(\mathcal{L}-\lambda)(\tau_{t_{0},x_{0}}\theta(t,x))-\mu(t,x)\alpha+\lambda \alpha\\
&>&\gamma\epsilon-\frac{\gamma\epsilon}{2}-\|\mu-\lambda\|_{\infty}\alpha\\
&>&0\\
\end{array} \]
for all $(t,x)\in B_{R}(x_{0})$. Thus, the strong maximum
principle yields that $\tilde{z}(t,x)= \alpha$ for all $t>t_{0}$ and $x\in\R^{N}$,
which contradicts $(\mathcal{L}-\lambda)(\tilde{z}-\alpha)>0$. $\Box$

\bigskip

\noindent {\bf Proof of Proposition \ref{cexlambda1'=0}.}
Assume that $u$ is a bounded positive continuous and entire solution of $(\ref{eqprinc})$. Then as $f(t,x,s)\leq \eta(t,x)s$ for all $(t,x,s)\in\R\times\R^{N}\times\R^{+}$, one has:
\[\partial_{t}u-\nabla\cdot (A\nabla u)+q\cdot \nabla u = f(t,x,u)\leq\eta (t,x)u.\]
As $u$ is positive and bounded, one can use $u$ as a test function in $(\ref{eigenlemmabounded})$. This gives $\lambda_{1}'(A,q,\eta)\leq 0$, which is a contradiction. $\Box$


\subsection{Existence of a minimal speed}

We now investigate the set
\[\mathcal{C}=\{c\in\R, \ \hbox{there exists some pulsating traveling front of speed} \ c\}.\]
In order to end the proof of Theorems $\ref{thmexistence}$ and $\ref{thmexistencepos}$ and to prove the existence of a minimal speed, it is only left to prove Proposition $\ref{vitferme}$, which yields that $\mathcal{C}$ is closed.

\bigskip

\noindent {\bf Proof of Proposition \ref{vitferme}.}
Consider a sequence $c_{n}\in\mathcal{C}$ which converges to some speed $c_{\infty}$. For all $n$, there exists a profile $\phi_{n}$ that satisfies equation $(\ref{eqtf})$ associated with the speed $c_{n}$. Up to some translation in $z$, one can assume that for all $n$,
\[\int_{(0,1)\times(0,T)\times C}\phi_n(z,t,x)dzdtdx=\frac{\min_{\R\times\R^{N}} p}{2}.\]
Proposition $\ref{removeregpos}$ yields that the sequence $(\phi_{n})$ is uniformly bounded in $W^{1,1}_{loc}(\R\times\R\times\R^{N})$. Thus one can assume, up to extraction, that this sequence converges almost everywhere and in $L^{1}_{loc}(\R\times\R\times\R^{N})$ to a function $\phi$. This function is a weak solution of equation $(\ref{eqtf})$ with $c=c_{\infty}$.
Using a similar method as in the proof of  Proposition $\ref{asymptotics2}$, we can prove that the good asymptotic behaviours hold when $z \rightarrow \pm\infty$. $\Box$

\bigskip

\noindent {\bf End of the proofs of Theorems \ref{thmexistence} and \ref{thmexistencepos}.} Under the hypotheses of Theorems $\ref{thmexistence}$ or $\ref{thmexistencepos}$, Proposition $\ref{vitferme}$ can be applied and thus the set $\mathcal{C}$ is closed. We also know that $\mathcal{C}$ contains the half-line $(c^{*}_{e}(A,q,\eta),+\infty)$. This gives the existence of a pulsating traveling front of speed $c=c^{*}_{e}(A,q,\mu)$ in Theorem $\ref{thmexistence}$ and of a pulsating traveling front of speed $c=c^{*}_{e}(A,q,\eta)$ in Theorem $\ref{thmexistencepos}$.$\Box$

\bigskip

We also easily get Theorem $\ref{mainresult}$ from this proposition:

\noindent {\bf Proof of Theorem \ref{mainresult}.}

We set $c^{*}_{e}=\inf \mathcal{C}$.
Theorem $\ref{thmexistencepos}$ yields that this set is not empty and contains the half-line $(c^{*}_{e}(A,q,\eta),+\infty)$. Proposition $\ref{thmnonexistence}$, which will be proved in section \ref{spreadingsection}, gives that it is bounded from below by $c^{*}_{e}(A,q,\mu)$. Thus the infimum is well-defined and
\[c^{*}_{e}(A,q,\mu)\leq c^{*}_{e} \leq c^{*}_{e}(A,q,\eta).\]

Proposition $\ref{vitferme}$ yields that $\mathcal{C}$ is closed and thus this infimum is in fact a $minimum$, that is, there exists a pulsating traveling front of speed $c^{*}_{e}$. This ends the proof.$\Box$


\section{Regularity of the pulsating traveling fronts}

This section is devoted to the proof of Theorem $\ref{mainresultcont}$. We assume that $f$ satisfies the hypotheses of Theorem $\ref{mainresult}$, in particular that $s\mapsto f(t,x,s)/s$ is nonincreasing, and we prove some uniform estimates in $W^{1,\infty}$ which guarantee that the profile $\phi$ we construct is Lipschitz continuous in $z$. We begin with the following lemma, which is of independent interest and which is true even if $f$ is not of KPP type:
\begin{lem} \label{phieps-infty}
Assume that $c>c^*_{\epsilon}(\mu)$ and consider $\phi_\epsilon$ a solution of (\ref{reginfinite}). Then 
\[\limsup_{z\rightarrow-\infty, (t,x)\in\R\times\R^{N}} \frac{\partial_{z}\phi_{\epsilon}}{\phi_{\epsilon}}= \lambda_{\epsilon}^{c} \ \hbox{or} \ \Lambda_{\epsilon}^{c}.\]
\end{lem}

\noindent {\bf Proof.}
Set $v=\displaystyle \frac{\partial_{z}\phi_{\epsilon}}{\phi_{\epsilon}}$. The Harnack inequality yields that $v$ is a bounded function. 
Furthermore, it satisfies:
\begin{equation}\begin{array}{rcl}
\partial_{t}v &=& \displaystyle \frac{\partial_{tz}\phi_{\epsilon}}{\phi_{\epsilon}}-
\frac{\partial_{t}\phi_{\epsilon}}{\phi_{\epsilon}}v\\
&&\\
\nabla v &=&\displaystyle  \frac{\partial_{z}\nabla\phi_{\epsilon}}{\phi_{\epsilon}}-
\frac{\nabla\phi_{\epsilon}}{\phi_{\epsilon}}v\\
&&\\
\partial_{z}v &=& \displaystyle \frac{\partial_{zz}\phi_{\epsilon}}{\phi_{\epsilon}}-
\frac{\partial_{z}\phi_{\epsilon}}{\phi_{\epsilon}}v\\
&&\\
\nabla\cdot (A\nabla v)&=&\displaystyle \frac{\partial_{z} \nabla\cdot
(A\nabla\phi_{\epsilon})}{\phi_{\epsilon}}-\frac{\nabla\cdot (A\nabla\phi_{\epsilon})}{\phi_{\epsilon}}v
-2\frac{\nabla\phi_{\epsilon}}{\phi_{\epsilon}}A\nabla v\\
&&\\
\partial_{zz}v&=&\displaystyle \frac{\partial_{zzz} \phi_{\epsilon}}{\phi_{\epsilon}}-\frac{\partial_{zz}\phi_{\epsilon}}{\phi_{\epsilon}}v
-2\frac{\partial_{z}\phi_{\epsilon}}{\phi_{\epsilon}}\partial_{z} v.\\
\end{array}\end{equation}

This computations yield that: 
\begin{equation}\label{eqfaiblez}\begin{array}{l}
\partial_{t}v-\nabla \cdot(A\nabla v)-2eA\nabla\partial_{z}v
 -eAe(1+\epsilon)\partial_{zz}v+q\cdot\nabla v+c\partial_{z}v\\
 \\
 +2\displaystyle \frac{\nabla\phi_{\epsilon}}{\phi_{\epsilon}}A\nabla v +2(1+\epsilon)\frac{\partial_{z}\phi_{\epsilon}}{\phi_{\epsilon}}Ae\partial_{z} v +2\frac{\nabla\phi_{\epsilon}}{\phi_{\epsilon}}Ae\partial_{z} v +2\frac{\partial_{z}\phi_{\epsilon}}{\phi_{\epsilon}}eA\nabla
 v\\
 \\
 =\displaystyle \frac{\partial_{z}(f(t,x,\phi_{\epsilon}))}{\phi_{\epsilon}}-\frac{f(t,x,\phi_{\epsilon})}{\phi_{\epsilon}}v\\
 \\
 =(f_u'(t,x,\phi_{\epsilon})-\displaystyle \frac{f(t,x,\phi_{\epsilon})}{\phi_{\epsilon}})v\leq 0.\\
 \end{array} \end{equation}

Next, set $m=\limsup_{z\rightarrow-\infty, (t,x)\in\R\times\R^{N}}\displaystyle \frac{\partial_{z}\phi_{\epsilon}}{\phi_{\epsilon}}$ and consider a sequence $(z_{n},t_{n},x_{n})$ such that $v(z_{n},t_{n},x_{n})\rightarrow m$ and $z_{n}\rightarrow-\infty$. For all $n$, there exists some $\overline{t}_{n}\in T\Z$ and $\overline{x}_{n}\in \Pi_{i=1}^{N}L_{i}\Z$ such that $s_{n}=t_{n}-\overline{t}_{n}\in [0,T]$ and $y_{n}=x_{n}-\overline{x}_{n}\in \overline{C}$. Up to extraction, we assume that $s_{n}\rightarrow s_{\infty}$ and $y_{n}\rightarrow y_{\infty}$. Set $\tilde{\phi}_{n}(t,x)=\displaystyle \frac{\phi_{\epsilon}(z+z_{n},t,x)}{\phi_{\epsilon}(z_{n},t_{n},x_{n})}$. This function satisfies:
\begin{equation} L_{\epsilon}\tilde{\phi}_{n}= \displaystyle \frac{1}{\phi_{\epsilon}(z_{n},t_{n},x_{n})}f(t,x,\tilde{\phi}_{n}\phi_{\epsilon}(z_{n},t_{n},x_{n})).
\end{equation}
The Schauder estimates yield that one may assume, up to extraction, that $\tilde{\phi}_{n}$ converges to some function $\tilde{\phi}_{\infty}$ in $\mathcal{C}^{2,1,2}_{loc}(\R\times\R\times\R^{N})$ as $n\rightarrow+\infty$, which is a solution of the linear equation:
\begin{equation}\label{eqtildephi} L_{\epsilon}\tilde{\phi}_{\infty}= \mu(t,x)\tilde{\phi}_{\infty},
\end{equation}
where $\mu(t,x)=f_u'(t,x,0)$ since $\phi_{\epsilon}(z_{n},t_{n},x_{n})\rightarrow 0$ as $n\rightarrow+\infty$. As $\tilde{\phi}_{\infty}(0,s_{\infty},y_{\infty})=1$ and $\tilde{\phi}_{\infty}$ is nonnegative, the strong maximum principle and the periodicity yield that $\tilde{\phi}_{\infty}$ is positive. 

Next, define $v_{n}(z,t,x)=v(z+z_{n},t,x)$, this function satisfies equation $(\ref{eqfaiblez})$, where $\phi_{\epsilon}(z,t,x)$ is replaced by $\phi_{\epsilon}(z+z_{n},t,x)$. Thus the Schauder estimates yield that the sequence $(v_{n})_{n}$ converges, up to extraction, to a function $v_{\infty}$ that satisfies:
\begin{equation}\begin{array}{l}
\partial_{t}v-\nabla \cdot(A\nabla v)-2eA\nabla\partial_{z}v
 -eAe(1+\epsilon)\partial_{zz}v+q\cdot\nabla v+c\partial_{z}v\\
 \\
 +2\displaystyle \frac{\nabla\phi_{\epsilon}}{\phi_{\epsilon}}A\nabla v +2(1+\epsilon)\frac{\partial_{z}\phi_{\epsilon}}{\phi_{\epsilon}}Ae\partial_{z} v +2\frac{\nabla\phi_{\epsilon}}{\phi_{\epsilon}}Ae\partial_{z} v +2\frac{\partial_{z}\phi_{\epsilon}}{\phi_{\epsilon}}eA\nabla
 v=0\\
 \end{array} \end{equation}

Furthermore, we know from the definition of $m$ that $v_{\infty}\geq m$ and that $v_{\infty}(0,s_{\infty},y_{\infty})=m$. Thus the strong parabolic maximum principle and the periodicity yield that $v_{\infty}\equiv m$ on $\R\times\R\times\R^{N}$.

As $\displaystyle \frac{\partial_{z}\tilde{\phi}_{n}}{\tilde{\phi}_{n}}=v_{n}$ for all $n$, one has $\displaystyle \frac{\partial_{z}\tilde{\phi}_{\infty}}{\tilde{\phi}_{\infty}}\equiv m$ and thus $\tilde{\phi}_{\infty}$ can be written $\tilde{\phi}_{\infty}(z,t,x)=\varphi (t,x) e^{mz}$, where $\varphi$ is periodic in $t$ and $x$ and positive. Reporting this in $(\ref{eqtildephi})$, one gets:
\[L_{m}\phi +mc\varphi -\epsilon m^{2}\varphi  = 0.\]
Thus $\varphi$ is some space-time periodic principal eigenfunction associated with $k_{m}$ and $k_{m}+mc-\epsilon m^{2}=0$. Thus we proved in Proposition $\ref{deflambdac}$ that this implies $m=\lambda_{c}^{\epsilon}(\mu)$ or $m=\Lambda_{c}^{\epsilon}(\mu)$. $\Box$

\bigskip

\begin{prop} \label{est-kpp-cont}
If $s\in\R^+\mapsto f(t,x,s)/s$ is nonincreasing for all $(t,x)\in\R\times\R^N$ and $c>c^*(\mu)$, then the derivatives $(\partial_{z}\phi_{\epsilon})_{\epsilon>0}$ are uniformly bounded
in $L^{\infty}(\R\times\R\times\R^{N})$.
\end{prop}

\noindent {\bf Remark.} This proof is not available if $s\mapsto f(t,x,s)/s$ is not nondecreasing. Actually, we need a sign for the zero order term of equation (\ref{eqfaiblez}). For example, taking $c$ large does not help. 

\bigskip

\noindent {\bf Proof.}
We now from the proof of the previous lemma that the function $v=\displaystyle \frac{\partial_{z}\phi_{\epsilon}}{\phi_{\epsilon}}$ satisfies $(\ref{eqfaiblez})$. As $f_u'(t,x,\phi_{\epsilon}(z,t,x))\phi\leq f(t,x,\phi_{\epsilon}(z,t,x))$ for all $(z,t,x)\in\R\times\R\times\R^{N}$ since $s\mapsto f(t,x,s)/s$ is nonincreasing, the weak maximum principle and the
periodicity yield that:
\begin{equation} \label{maxprincv}0\leq
v(z,t,x) \leq\max_{(t,x)\in\R\times\R^{N}}
\{v(-a,t,x),v(a,t,x)\}.\end{equation} 

As $\phi_{\epsilon}$ is nondecreasing to $p$ as $z\rightarrow +\infty$ and as it satisfies a parabolic equation, the Schauder estimates yield that $\partial_{z}\phi_{\epsilon}\rightarrow 0$ as $z\rightarrow+\infty$. Thus $v(a,t,x)\rightarrow 0$ as $a\rightarrow+\infty$ uniformly with respect to $(t,x)\in\R\times\R^{N}$.
Furthermore, Lemma $\ref{phieps-infty}$ gives that $\limsup_{z\rightarrow+\infty, (t,x)\in\R\times\R^{N}} v(-a,t,x) \leq\Lambda_{c}^{\epsilon}$ and this quantity is uniformly bounded by some constant $R_{c}$ which does not depend on $\epsilon$ from Lemma $\ref{contlambdaeps}$.

Thus the right-hand
side of $(\ref{maxprincv})$ is uniformly bounded with respect to $\epsilon$ and $a$ by a positive
constant $R_{c}$. Finally, for all $(z,t,x)\in\R\times\R\times\R^{N}$,
we have:
\[0\leq \partial_{z}\phi_{\epsilon}(z,t,x) \leq R_{c}\|p\|_{\infty}.\]
$\Box$

\bigskip

We are now able to prove Theorem $\ref{mainresultcont}$.

\smallskip

\noindent {\bf Proof of Theorem \ref{mainresultcont}.}
It is only left to prove that one can modify the proof of Theorem $\ref{thmexistence}$ in order to get a Lipschitz continuous profile $\phi$. First, assume that $c>c^*(\mu)$. 
Fix a compact set $K\subset \R$ and let
\[\begin{array}{rcccc}
\Phi_\epsilon &:& K &\rightarrow & L^2_{per}(\R\times\R^N)\\
&& z&\mapsto& \big( (t,x)\mapsto \phi_\epsilon (z,t,x)\big)\\
\end{array} \]
where $L^2_{per}(\R\times\R^N)$ is the space of the functions that are space-time periodic and that belong to $L^2_{loc}(\R\times\R^N)$. Propositions $\ref{est-kpp-cont}$ and $\ref{removeregpos}$ yield that for all $z\in K$, the family $(\Phi_\epsilon(z))_{\epsilon>0}$ is uniformly bounded in $H^1_{per}(\R\times\R^N)$. Hence $(\Phi_\epsilon(z))_{\epsilon>0}$ is relatively compact in $L^2_{per}(\R\times\R^N)$ for all $z\in K$. 

Moreover, Proposition $\ref{est-kpp-cont}$ yields that the family $(\Phi_\epsilon)_{\epsilon>0}$ is equicontinuous. The Ascoli theorem gives that it is a relatively compact family and thus we can assume that it converges to some $\Phi\in \mathcal{C}^0\big(K,L^2_{per}(\R\times\R^N)\big)$ as $\epsilon\rightarrow 0$. Using a diagonal extraction process , we can assume that $\phi_\epsilon\rightarrow \phi$ in $\mathcal{C}^{0}_{loc}(\R, L^{2}_{per}(\R\times\R^{N}))$ as $\epsilon\rightarrow 0$, where $\phi(z,\cdot,\cdot)=\Phi(z)$. As the estimate given by Proposition $\ref{est-kpp-cont}$ is uniform in $z\in\R$, $\phi$ belongs to $W^{1,\infty}(\R, L^{2}_{per}(\R\times\R^{N}))$.

Setting $u(y,t,x)=\phi(y+x\cdot e+ct,t,x)\in W^{1,\infty}(\R, L^{2}_{per}(\R\times\R^{N}))$, one gets a parametrized family of functions $u_{y}: (t,x)\mapsto u(y,t,x)$ such that for all $y$, $u_{y}$ satisfies $(\ref{eqprinc})$ in the sense of distribution since $\phi$ satisfies $(\ref{eqtf})$ and $y\mapsto u_{y}$ is continuous. Thus for all $y$, $u_{y}\in\mathcal{C}^{1,2}(\R\times\R^{N})$ from the Schauder parabolic estimates. Thus $\phi$ is Lipschitz continuous with respect to $(z,t,x)\in\R\times\R\times\R^N$. 

If $c=c^*(\mu)$, we consider a sequence $(c_n)_n$ as in the proof of Proposition \ref{vitferme}. As the estimates of Proposition \ref{est-kpp-cont} are uniform with respect to the sequence $(c_n)_n$, the sequence of profiles associated with the speeds $(c_n)_n$ is uniformly bounded in appropriate norms and one can pass to the limit as previously. This gives a Lipschitz continuous pulsating traveling front of speed $c^*(\mu)$. $\Box$


\section{Spreading properties}\label{spreadingsection}

We now prove the spreading properties for front-like initial data. The aim of this section is also to get the nonexistence Theorem $\ref{thmnonexistence}$.

\begin{lem} \label{eigenvaluecomplex}
There exists $c'<c^{*}_{e}(A,q,\mu)$ such that for all $c\in(c',c^{*}(\mu))$, there exists a complex $\lambda \in \mathbb{C}\backslash\R$ and a solution $\psi \in \mathcal{C}^{1,2}(\R,\R^{N})$ of:
\begin{equation}  \left\{ \begin{array}{l}
\partial_{t}\psi - \nabla\cdot(A\nabla \psi) -2\lambda e A \nabla \psi+q\cdot\nabla\psi-(\lambda^{2} e A e-\lambda c+\lambda\nabla\cdot(A e)+\mu-\lambda q\cdot e)\psi=0,\\
\psi \ \hbox{is periodic in $(t,x)$},\\
Re(\psi)>0.\\
\end{array} \right. \end{equation}
\end{lem}

\noindent {\bf Proof.}
Set $\lambda^{*}=\lambda_{c^{*}_{e}(A,q,\mu)}$. The family of operators $L_{\lambda}$ depends analytically on $\lambda$, in the sense of Kato. From the Kato-Rellich theorem \cite{Kato}, there exists a neighborhood $V$ of $\lambda^{*}$ in $\mathbb{C}$, such that there exists a simple eigenvalue $\tilde{k}_{\lambda}(\mu)$ continuing $k_{\lambda}(\mu)$ on all $V$ analytically and a family of eigenfunctions $\psi_{\lambda}$ analytic in $\lambda$, where $\psi_{\lambda^{*}}$ is the positive principal eigenfunction associated with $c^{*}_{e}(A,q,\mu)$.

Set $F_{c}(\lambda)=\tilde{k}_{\lambda}(\mu)+\lambda c$. This function is analytic in $\lambda$ and converges locally uniformly to $F_{c^{*}_{e}(A,q,\mu)}$ as $c\rightarrow c^{*}_{e}(A,q,\mu)$. As $F_{c^{*}_{e}(A,q,\mu)}(\lambda^{*})=0$, the Rouch\'e theorem yields that there exists some neighborhood $V$ of $c^{*}_{e}(A,q,\mu)$ such that for all $c \in V$, there exists some $\lambda_{c}\in\mathbb{C}$ such that $F_{c}(\lambda_{c})=0$ and $\lambda_{c}\rightarrow \lambda^{*}$ as  $c\rightarrow c^{*}_{e}(A,q,\mu)$. 

Using the classical Schauder estimates, one can prove that $\psi_{\lambda_{c}}\rightarrow \psi_{\lambda^{*}}$ uniformly in $t$ and $x$. Thus $Re(\psi_{\lambda_{c}})\rightarrow \psi_{\lambda^{*}}>0$ and taking $V$ small enough, we can assume that $Re(\psi_{\lambda_{c}})>0$ for all $c \in V$.
Lastly, if $c<c^{*}_{e}(A,q,\mu)$, it is impossible to have $\lambda_{c}\in \R$. Otherwise, this would contradict the definition of $c^{*}_{e}(A,q,\mu)$. This ends the proof of the lemma.$\Box$

\bigskip

\noindent {\bf Proof of Proposition \ref{spreadingsub}.}
First, we assume that $c'<c<c^{*}_{e}(A,q,\mu)$. We know that $c^{*}(\mu-\delta)\rightarrow c^{*}(\mu)$ as $\delta \rightarrow 0$, so that one can fix $\delta >0$ such that $c<c^{*}(\mu-\delta)<c^{*}(\mu)$. As $f$ is of class $\mathcal{C}^{1}$ in $\R\times\R^{N}\times[0,\beta]$ for a given positive $\beta$, there exists a positive constant $\epsilon$ such that:
\[ \forall (t,x,s) \in \R\times\R^{N}\times [0,\epsilon], f(t,x,s)\geq (\mu(t,x)-\delta)s.\]

We set $\psi$ associated with $c$ given by Lemma $\ref{eigenvaluecomplex}$ and we consider the function:
\begin{equation} w_{0}(t,x) = Re(e^{\lambda(x\cdot e+ct)}\psi(t,x)).\end{equation}
Next, one has:
\begin{equation} \label{reecriturew} w_{0}(t,x)=e^{\lambda_{r}(x\cdot e+ct)}[ \psi_{r}cos(\lambda_{i}(x\cdot e+ct))+\psi_{i}sin(\lambda_{i}(x\cdot e+ct)) ] \end{equation}
where $\psi_{i},\psi_{r},\lambda_{i},\lambda_{r}$ denote the imaginary and real parts of $\lambda$ and $\psi$. 
For all $n \in \Z$, if $(e\cdot x+ ct) = 2n \pi/\lambda_{i}$, then $w_{0}(t,x)>0$. Similarly, for all $n \in \Z$, if $(e\cdot x+ ct) = (2n+1) \pi/\lambda_{i}$, then $w_{0}(t,x)<0$. Thus, it follows from $(\ref{reecriturew})$ that there exists an interval $[b_{1},b_{2}]\subset \R$ and an unbounded domain $D\subset \R\times\R^{N}$ such that:
\begin{equation} \left\{ \begin{array}{l}
D \subset \big\{ (t,x)\in \R\times \R^{N}, x\cdot e+ct \in [b_{1},b_{2}]\big\},\\
0<w_{0}(t,x)<\epsilon, \hbox{for all} \ (t,x)\in D,\\
w_{0}(t,x)=0, \ \hbox{for} \ (t,x)\in\partial D.\\
\end{array} \right. \end{equation}

Set $w$ the function:
\begin{equation} 
w(t,x)=\left\{ \begin{array}{l}
  w_{0}(t,x) \ \hbox{if} \ (t,x)\in D,\\
  0 \ \hbox{otherwise}. \\
\end{array} \right. \end{equation}
This function satisfies the inequation:
\[\partial_{t}w-\nabla\cdot (A\nabla w)+q\cdot\nabla w = (\mu-\delta)w\leq f(t,x,w) \ \hbox{for all} \ (t,x)\in D.\]

Assume first that $u_{0}(x)=w(0,x)$. In this case the parabolic maximum principle yields that $u \geq w$. Set $v(t,x)=u(t,x-cte)$, $B(t,x)=A(t,x-cte)$, $r(t,x)=q(t,x-cte)$ and $g(t,x,s)=f(t,x-cte,s)$. The function $v$ is the solution of 
\[\partial_t v -\nabla\cdot (B(t,x)\nabla v) -r(t,x) \cdot \nabla v = g(t,x,v) \hbox{ in } \R\times\R^N.\]
Moreover, $B$, $r$ and $g$ are almost periodic in $(t,x)\in\R\times\R^N$. That is, for any sequence $(t_n,x_n)$ in $\R\times\R^N$, there exists a subsequence $(t_{n'},x_{n'})$ such that the sequences $(B(t+t_{n'},x+x_{n'}))_{n'}$, $(r(t+t_{n'},x+x_{n'}))_{n'}$ and $(g(t+t_{n'},x+x_{n'}))_{n'}$ converge uniformly in $(t,x)\in\R\times\R^N$ and locally uniformly in $s\geq 0$. 

Take an arbitrary sequence $t_{n}\rightarrow +\infty$ as $n\rightarrow +\infty$. Set $v_{n}(t,x)=v(t+t_{n},x)$, this function satisfies:
\[\partial_{t}v_{n}-\nabla\cdot (B(t+s_{n},x)\nabla v_{n}) +r(t+s_{n},x)\cdot\nabla v_{n}-c e\cdot\nabla v_{n}=g(t+s_{n},x,v_{n}) \ \hbox{in} [-t_{n},+\infty)\times\R^{N}.\]

Up to extraction, one may assume that $(B(t+t_{n},x))_{n}$, $(r(t+t_{n},x))_{n}$ and $(g(t+t_{n},x))_{n}$ converge uniformly in $(t,x)\in\R\times\R^N$ and locally uniformly in $s\geq 0$ to some $B_\infty, r_\infty$ and $g_\infty$ as $n\rightarrow+\infty$. The classical Schauder estimates then yield that one may find a subsequence $(v_{n'})$ that uniformly converges on any compact subset to a function $v_{\infty}$ in $\mathcal{C}^{1,2}_{loc}(\R\times\R^{N})$. The function $v_{\infty}$ is nonnegative and satisfies:
\[\partial_{t}v_{\infty}-\nabla\cdot (B_\infty(t,x)\nabla v_{\infty})+r_\infty(t,x)\cdot\nabla v_{\infty}-c e\cdot\nabla v_{\infty}=g_\infty(t,x,v_{\infty}) \ \hbox{in} \ \R\times\R^{N}.\]

Furthermore, for all $n$, one has:
\[ \begin{array}{l}
v_{n}(t,x)=u(t+t_{n},x-c(t+t_{n})e)\geq w(t+t_{n},x-c(t+t_{n})e)\\
\geq e^{\lambda_{r}(x\cdot e)}[ \psi_{r}(t+t_{n},x-c(t+t_{n})e)cos(\lambda_{i}(x\cdot e))+\psi_{i}(t+t_{n},x-c(t+t_{n})e)sin(\lambda_{i}(x\cdot e))].\\ \end{array} \]
Thus, taking $x_{0}=\frac{2n\pi}{\lambda_{i}}e$ and using the positivity and the periodicity of $\psi_{r}$, one gets $\inf_{n\in\mathbb{N}}\inf_{t\in\R} v_n(t,x_0)>0$, which yields that $\inf_{t\in\R}v_\infty(t,x_0)>0$.
The Krylov-Safonov-Harnack inequality yields that $\inf_{(t,x)\in\R\times \overline{C}}v_{\infty}(t,x)>0$. As $v_{n}$ is periodic in $x$ for all $n$, the function $v_{\infty}$ is also periodic in $x$ and then $\inf_{(t,x)\in\R\times \R^{N}}v_{\infty}(t,x)>0$. Hypothesis $\ref{uniquenesshyp}$ then yields that $v_{\infty}(t,x)=\lim_{n'\rightarrow+\infty} p(t+t_{n'},x-c(t+t_{n'})e)$ and thus
\[v_{n'}(t,x)-p(t+t_{n'},x-c(t+t_{n'})e)\rightarrow 0\]
as $n'\rightarrow +\infty$, uniformly on every compact subset.
Finally, the classical procedure yields that:
\begin{equation} \label{convergencespreadingsub} u(t,x-cte)-p(t,x-cte) \rightarrow 0 \ \hbox{as} \ t\rightarrow +\infty .\end{equation}
uniformly on any compact subset.

To sum up, we have constructed an initial datum $w(0,.)$ with compact support such that the solution $u$ associated with this initial datum satisfies $(\ref{convergencespreadingsub})$. Furthermore, multiplying the function $\psi$ by a positive constant, one can take an arbitrary small supremum norm for $w(0,\cdot)$. Applying the maximum principle principle, we generalize this result to any initial datum $u_{0}$ such that there exists $(a_{1},a_{2})\in\R^{2}$ such that $\inf_{x\cdot e\in [a_{1},a_{2}]}u_{0}(x) >0$.$\Box$

\bigskip

Next, we prove the result for any speed $-c^{*}_{-e}(A,q,\mu)<c<c^{*}_{e}(A,q,\mu)$. 
Set:
\[\Omega=\{(t,x)\in\R^{+}\times\R^{N}, a_{1}-c^{*}_{-e}(A,q,\mu)t\leq x\cdot e \leq a_{2}+c^{*}_{e}(A,q,\mu)t\}\]
One has $\inf_{a_{1}\leq x\cdot e\leq a_{2}}u_{0}(x)>0$. The previous case yields that $u(t,x-cte)-p(t,x-cte)\rightarrow 0$ as $t\rightarrow +\infty$ when $c$ is close to $c^{*}_{e}(A,q,\mu)$ and $u(t,x+cte)-p(t,x+cte)\rightarrow 0$ as $t\rightarrow +\infty$ when $c$ is close to $-c^{*}_{-e}(A,q,\mu)$, where $p$ is positive and periodic in $t$ and $x$. Thus there exists some $\epsilon > 0$ such that for all $(t,x)\in\partial\Omega, u(t,x)>\epsilon$. We need a modified maximum principle in order to get an estimate in the whole set $\Omega$. As $\Omega$ is not a cylinder, we cannot apply the classical weak maximum principle. In fact, it is possible to extend this maximum principle to the cone $\Omega$ and there is no particular issue but, by sake of completeness, we prove that this extension works well here:

\begin{lem} 
Assume that $z$ satisfies:
\[ \left\{ \begin{array}{r}\partial_{t}z-\nabla\cdot (A\nabla z)+q\cdot \nabla z +b z \geq 0 \ \hbox{in} \ \Omega\\
z \geq 0 \ \hbox{in} \ \partial\Omega\\
\end{array}\right. \]
where $b$ is a bounded continuous function. then one has $z\geq 0$ in $\Omega$.
\end{lem}

\noindent {\bf Proof.}
Assume first that $b > 0$. Set $\Omega_{\tau}=\Omega \cap \left\{t\leq\tau\right\}$ and assume that there exists $(t,x)\in \overline{\Omega_{\tau}}$ such that $z(t,x)<0$. Take $(t_{0},x_{0}) \in \overline{\Omega_{\tau}}$ such that $z(t_{0},x_{0})=\min_{(t,x)\in\overline{\Omega_{\tau}}}z(t,x)<0$. One necessarily has $(t_{0},x_{0}) \in \Omega_{\tau}$ and thus:
\[\nabla z(t_{0},x_{0})=0, \ \nabla\cdot (A\nabla z)(t_{0},x_{0}) \geq 0, \ b(t_{0},x_{0})z(t_{0},x_{0}) <0.\]
This leads to:
\[\partial_{t}z(t_{0},x_{0})> 0.\]
But the def intion of the minimum yields that for all $0\leq t\leq t_{0}$, if $(t,x_{0})\in \Omega$, one has $z(t,x_{0})\geq z(t_{0},x_{0})$. As $t_{0}>0$, for $\epsilon$ small enough, one has $(t_{0}-\epsilon,x_{0}) \in \Omega$. Thus it is possible to differentiate the inequality, which gives $\partial_{t}z(t_{0},x_{0})\leq 0$. This is a contradiction. Thus for all $\tau>0$, one has $\min_{\Omega_{\tau}}z \geq 0$ and then $z\geq 0$ in $\Omega$.

If $b$ is not positive, set $z_{1}(t,x)=e^{-(\|b\|_{\infty}+1)t}z(t,x)$ for all $(t,x)\in \Omega$. This function satisfies:
\[\partial_{t}z_{1}-\nabla\cdot (A\nabla z_{1})+q\cdot z_{1} +(b+\|b\|_{\infty}+1) z_{1}=(\partial_{t}z-\nabla\cdot (A\nabla z)+q\cdot z +b z)e^{-(\|b\|_{\infty}+1)t}\geq 0,\]
and for all $(t,x)\in\partial\Omega$, one has $z_{1}(t,x) \geq 0$. As $b+\|b\|_{\infty}+1 >0$, the first case yields that $z_{1}\geq 0$ and then $z\geq 0$. This ends the proof.$\Box$

\bigskip

In order to apply this lemma, take $\psi=\psi_{0}$ a periodic principal eigenfunction associated with $L_0$ such that $\|\psi\|_{\infty}<\epsilon $. Set $z=u-\psi$ and $b(t,x)=\frac{f(t,x,u)-f(t,x,\psi)}{u-\psi}$. As $f$ is Lipschitz continuous in $u$ uniformly in $(t,x)$, the function $b$ is bounded. The function $z$ satisfies the equation:
\[\partial_{t}z-\nabla\cdot (A\nabla z)+q\cdot z -b z=0.\]
Thus, the hypothesis of the previous lemma are satisfied and one has $z\geq 0$, that is, $u\geq\psi$ in $\Omega$.

Take now $c\in(-c^{*}_{-e}(A,q,\mu),c^{*}_{e}(A,q,\mu))$, as $(t,x-cte) \in \Omega$, one has 
\[\inf_{t\in\R^{+}, x\in\R^{N}}u(t,x-cte) \geq\epsilon.\] 
Take $t_{n}\rightarrow\infty$ and set $v_{n}(t,x)=u(t+t_{n},x-c(t+t_{n})e)$, up to extraction, one may assume that $v_{n}$ converge to a function $v_{\infty}$ in $\mathcal{C}^{1,2}_{loc}(\R\times\R^{N})$. The function $v_{\infty}$ is an entire bounded solution of an equation of type $(\ref{eqprinc})$ and satisfies $\inf_{\R\times\R^{N}}v_{\infty} \geq \epsilon>0$. Furthermore, it is periodic in $x$. Hypothesis $\ref{uniquenesshyp}$ yields that $v_{\infty}\equiv p$. The classical extraction arguments concludes the proof.$\Box$

\bigskip

\noindent {\bf Proof of Proposition \ref{thmnonexistence}.}
Assume that such a pulsating traveling front $u$ of speed $c<c^*(\mu)$ does exists and set $\phi$ its profile. Up to some shift of $\phi$ in $z$, we can assume that $u$ satisfies (\ref{eqprinc}). Then as $\phi(x\cdot e+ct,t,x)-p(t,x) \rightarrow 0$ uniformly in $x$ as $t\rightarrow +\infty$ and $p$ is a positive periodic function, $u(t,x)=\phi(x\cdot e+ct,t,x)$ satisfies the hypothesis of Proposition $\ref{spreadingsub}$. Thus, taking $c' \in (c,c^{*}_{e}(A,q,\mu))$ such that $c'\geq c^{**}(\mu)$, one gets:
\[u(t,x-c'te)-p(t,x)=\phi(x\cdot e-(c'-c)te,t,x-c'te)-p(t,x)\rightarrow 0 \ \hbox{as} \ t\rightarrow +\infty \]
In the other hand, as $c'-c>0$, one has $\phi(x\cdot e-(c'-c)te,t,x-c'te)\rightarrow 0$ as $t\rightarrow +\infty$ uniformly in $x$. As $p$ is positive, this gives a contradiction.$\Box$

\bigskip

\noindent {\bf Proof of Proposition \ref{spreadingsup}.}
Take $u_{0}$ an initial datum that satisfies the hypotheses and $c'\in(c^{*}_{e}(A,q,\eta),c)$. Set $v(t,x)=\psi_{\lambda_{c'}(\eta)}(t,x)e^{\lambda_{c'}(\eta)(x\cdot e+c't)}$, where $\psi_{\lambda_{c'}(\eta)}$ is the periodic principal eigenvalue normalized by $\|\psi_{\lambda_{c'}}\|_{\infty}=1$. Since $c'>c^{*}(\eta)$, one has $\lambda_{c'}(\eta)<\lambda^{c^{*}(\eta)}(\eta)$ and then the hypotheses yield that there exist two positive constants $A, C$ such that:
\[u_{0}(x)\leq C v(0,x) \ \hbox{if} \ x\cdot e<-A.\]
Thus, as one can increase $C$, for all $x\in \R^{N}$, one has $u_{0}(x)\leq C v(0,x)$. The function $Cv$ is a subsolution of equation $(\ref {eqprinc})$ and the maximum principle thus gives $u(t,x)\leq C v(t,x)$ for all $(t,x)$. 

Finally, one has:
\[u(t,x-cte)\leq C e^{\lambda_{c'}(\eta)(x\cdot e+(c'-c)t)} \rightarrow 0 \ \hbox{as} \ t\rightarrow +\infty\]
uniformly in $x\cdot e\leq -B$.$\Box$

\subsection*{Acknowledgments} 
It is my pleasure to thank professors F. Hamel and H. Berestycki for the attention they paid to this work. Professor F. Hamel suggested Lemma $\ref{thetalem}$ and the proof of Theorem $\ref{thmexistencepos}$, which considerably simplified all the proofs and weakened the hypotheses. I would also like to thank J.-M. Roquejoffre for his suggestions of improvements, professor J. Nolen for the precisions he gave me about Definition $\ref{deftravelingfrontsNolen}$, professor A. Zlatos for the suggestions he did after having read the preprint of this paper and the University of Chicago, where this work has started and ended.

\bibliographystyle{plain}
\bibliography{biblio2}

\begin{thebibliography}{10}

\bibitem{Alikakos}
N.~Alikakos, P.~W. Bates, and X.~Chen.
\newblock Traveling waves in a time periodic structure and a singular
  perturbation problem.
\newblock {\em Transactions of AMS}, 351:2777--2805, 1999.

\bibitem{Aronson}
D.G. Aronson and H.F. Weinberger.
\newblock Multidimensional nonlinear diffusions arising in population genetics.
\newblock {\em Adv. Math.}, 30:33--76, 1978.

\bibitem{Generalwaves2}
H.~Berestycki and F.~Hamel.
\newblock On a general definition of transition waves and their properties.
\newblock {\em preprint}.

\bibitem{Frontexcitable}
H.~Berestycki and F.~Hamel.
\newblock Front propagation in periodic excitable media.
\newblock {\em Comm. Pure Appl. Math.}, 55:949--1032, 2002.

\bibitem{Generalwaves}
H.~Berestycki and F.~Hamel.
\newblock Generalized travelling waves for reaction-diffusion equations.
\newblock {\em Perspectives in Nonlinear Partial Differential Equations. In
  honor of H. Brezis, Contemp. Math. 446, Amer. Math. Soc.}, pages 101--123,
  2007.

\bibitem{BHKR}
H.~Berestycki, F.~Hamel, A.~Kiselev, and L.~Ryzhik.
\newblock Quenching and propagation in kpp reaction-diffusion equations with a
  heat loss.
\newblock {\em Arch. Ration. Mech. Anal.}, 178:57--80, 2005.

\bibitem{Base1}
H.~Berestycki, F.~Hamel, and L.Roques.
\newblock Analysis of the periodically fragmented environment model : 1 -
  influence of periodic heterogeneous environment on species persistence.
\newblock {\em J. Math. Biol.}, 51:75--113, 2005.

\bibitem{Base2}
H.~Berestycki, F.~Hamel, and L.Roques.
\newblock Analysis of the periodically fragmented environment model : 2 -
  biological invasions and pulsating travelling fronts.
\newblock {\em J. Math. Pures Appl.}, 84:1101--1146, 2005.

\bibitem{posspreadspeed}
H.~Berestycki, F.~Hamel, and G.~Nadin.
\newblock Asymptotic spreading in diffusive excitable media.
\newblock {\em to appear in J. Func. Anal.}, 2007.

\bibitem{Larrouturou}
H.~Berestycki, B.~Larrouturou, and P.~L. Lions.
\newblock Multi-dimensional travelling-wave solutions of a flame propagation
  model.
\newblock {\em Arch. Rational Mech. Anal.}, 111(1):33--49, 1990.

\bibitem{Nirenberg}
H.~Berestycki and L.~Nirenberg.
\newblock Travelling fronts in cylinders.
\newblock {\em Ann. Inst. H. Poincar\'e}, 9:497--572, 1992.

\bibitem{TheseBages}
M.~Bag\` es.
\newblock \'equations de r\'eaction-diffusion de type kpp : ondes pulsatoires,
  dynamiques non-triviales et applications.
\newblock {\em PhD thesis}, 2007.

\bibitem{Fisher}
R.~A. Fisher.
\newblock The advance of advantageous genes.
\newblock {\em Ann. Eugenics}, 7:335--369, 1937.

\bibitem{Freidlin2}
M.~Freidlin.
\newblock On wave front propagation in periodic media.
\newblock {\em In: Stochastic analysis and applications, ed. M. Pinsky,
  Adavances in Probability and related topics}, 7:147--166, 1984.

\bibitem{Gartner}
M.~Freidlin and J.~Gartner.
\newblock On the propagation of concentration waves in periodic and random
  media.
\newblock {\em Sov. Math. Dokl.}, 20:1282--1286, 1979.

\bibitem{Frejacques}
G.~Frejacques.
\newblock Travelling waves in infinite cylinders with time-periodic
  coefficients.
\newblock {\em PhD Thesis}.

\bibitem{GuoHamel2}
J.-S. Guo and F.~Hamel.
\newblock Propagation and slowdown in hostile environments.
\newblock {\em in preparation}.

\bibitem{GuoHamel1}
J.-S. Guo and F.~Hamel.
\newblock Front propagation for discrete periodic monostable equations.
\newblock {\em Math. Ann.}, 335:489--525, 2006.

\bibitem{Hamelmonoticity}
F.~Hamel.
\newblock Qualitative properties of monostable pulsating fronts : exponential
  decay and monoticity.
\newblock {\em preprint}, 2007.

\bibitem{Hamelmonotonicity}
F.~Hamel.
\newblock Qualitative properties of monostable pulsating fronts : exponential
  decay and monoticity.
\newblock {\em J. Math. Pures Appl.}, 89:355--399, 2008.

\bibitem{HamelRyzhik}
F.~Hamel and L.~Ryzhik.
\newblock Non-adiabatic kpp fronts with an arbitrary lewis number.
\newblock {\em Nonlinearity}, 18:2881--2902, 2005.

\bibitem{Hormander}
L.~Hormander.
\newblock Hypoelliptic second order differential equations.
\newblock {\em Acta Math.}, 119:147171, 1967.

\bibitem{Kato}
T.~Kato.
\newblock Perturbation theory for linear operators.
\newblock {\em Springer Verlag, Berlin}, 1980.

\bibitem{Kohn}
J.~J. Kohn and L.~Nirenberg.
\newblock Degenerate elliptic-parabolic equations of second order.
\newblock {\em Comm. Pure Appl. Math.}, 20:797872, 1967.

\bibitem{KPP}
A.N. Kolmogorov, I.G. Petrovsky, and N.S. Piskunov.
\newblock Etude de l \'equation de la diffusion avec croissance de la
  quantit\'e de mati\`ere et son application \`a un probl\`eme biologique.
\newblock {\em Bulletin Universit\'e d'Etat \`a Moscou (Bjul. Moskowskogo Gos.
  Univ.)}, pages 1--26, 1937.

\bibitem{Roquejoffre}
J.-F. Mallordy and J.-M. Roquejoffre.
\newblock A parabolic equation of the kpp type in higher dimensions.
\newblock {\em SIAM J. Math. Anal.}, 26(1).

\bibitem{Matanodef}
H.~Matano.
\newblock Traveling waves in spatially inhomogeneous diffusive media.
\newblock {\em Oral communications}.

\bibitem{MelletRoquejoffre}
A.~Mellet and J.-M. Roquejoffre.
\newblock Construction dondes g\'en\'eralis\'ees pour le modele 1d scalaire \`a
  temp\'erature dignition.
\newblock {\em preprint}, 2007.

\bibitem{noteexistence}
G.~Nadin.
\newblock Reaction-diffusion equations in space-time periodic media.
\newblock {\em C. R. Acad. Sci. Paris, Ser. I}, 345(9).

\bibitem{existence}
G.~Nadin.
\newblock Existence and uniqueness of the solution of a space-time periodic
  reaction-diffusion equation.
\newblock {\em submitted}, 2007.

\bibitem{eigenvalue}
G.~Nadin.
\newblock The principal eigenvalue of a space-time periodic parabolic operator.
\newblock {\em to appear in Ann. Mat. Pura Appl.}, 2008.

\bibitem{notetf}
G.~Nadin.
\newblock Pulsating traveling fronts in space-time periodic media.
\newblock {\em C.R. Acad. Sci. Paris I}, 346:951--956, 2008.

\bibitem{dependance}
G.~Nadin.
\newblock Some remarks on the dependance relations between the speed of
  propagation for space-time periodic kpp equations and the coefficients.
\newblock {\em preprint}, 2008.

\bibitem{NadinRossi}
G.~Nadin and L.~Rossi.
\newblock Traveling fronts in space-time almost periodic media.
\newblock {\em in preparation}.

\bibitem{Nolen}
J.~Nolen, M.~Rudd, and J.~Xin.
\newblock Existence of kpp fronts in spatially-temporally periodic advection
  and variational principle for propagation speeds.
\newblock {\em Dynamics of PDE}, 2(1).

\bibitem{NolenRyzhik}
J.~Nolen and L.~Ryzhik.
\newblock Traveling waves in a one-dimensional random medium.
\newblock {\em preprint}, 2007.

\bibitem{Nolenshear}
J.~Nolen and J.~Xin.
\newblock Existence of kpp type fronts in space-time periodic shear flows and a
  study of minimal speeds based on variational principle.
\newblock {\em Discrete and Continuous Dynamical Systems}, 13(5).

\bibitem{NolenXin1d}
J.~Nolen and J.~Xin.
\newblock Kpp fronts in 1d random drift.
\newblock {\em to appear in Discrete and Continuous Dynamical Systems B}, 2008.

\bibitem{Shengeneral}
W.~Shen.
\newblock Traveling waves in time dependent bistable media.
\newblock {\em Diff. Int. Eq.}, 19(3).

\bibitem{Shen1}
W.~Shen.
\newblock Traveling waves in time almost periodic structures governed by
  bistable nonlinearities, 1. stability and uniqueness.
\newblock {\em J. Diff. Eq.}, 159:1--55, 1999.

\bibitem{Shen2}
W.~Shen.
\newblock Traveling waves in time almost periodic structures governed by
  bistable nonlinearities, 2. existence.
\newblock {\em J. Diff. Eq.}, 159:55--101, 1999.

\bibitem{Biologicalinvasions}
N.~Shigesada and K.~Kawasaki.
\newblock Biological invasions : theory and practice.
\newblock {\em Oxford Series in Ecology and Evolution, Oxford : Oxford
  University Press}, 1997.

\bibitem{Shigesada1}
N.~Shigesada, K.~Kawasaki, and E.~Teramoto.
\newblock Traveling periodic waves in heterogeneous environments.
\newblock {\em Theor. Population Biol.}, 30:143--160, 1986.

\bibitem{Weinberger}
H.~Weinberger.
\newblock On spreading speed and travelling waves for growth and migration
  models in a periodic habitat.
\newblock {\em J. Math. Biol.}, 45:511--548, 2002.

\bibitem{Xin}
J.~Xin.
\newblock Existence of planar flame fronts in convective-diffusive periodic
  media.
\newblock {\em Arch. Ration. Mech. Anal.}, 121:205--233, 1992.

\end{thebibliography}

\end{document}